\documentclass[11pt,a4paper]{article}
\usepackage{amsmath}
\usepackage{bbding}

\usepackage{amsmath,amsfonts,amssymb,amsthm}

\newcommand{\diam}{\operatorname{diam}}

\newcommand{\counte}{equation}%%%设置定理与公式公用计数器
\newtheorem{theorem}[\counte]{\bf Theorem}%[section]

\newtheorem{prop}[\counte]{\bf Proposition}
\newtheorem{lemma}[\counte]{\bf Lemma}

\newtheorem{coro}[\counte]{\bf Corollary}

\newtheorem{remark}[\counte]{\bf Remark}

\allowdisplaybreaks

\numberwithin{equation}{section}%%公式随section重新计数
\renewcommand{\thefootnote}{\fnsymbol{footnote}}

\begin{document}

\renewcommand{\thefootnote}{\arabic{footnote}}

\centerline{\bf\Large An Isometrical ${\Bbb C\Bbb
P}^{n}$-Theorem \footnote{Supported by NSFC 11471039. \hfill{$\,$}}}

\vskip5mm

\centerline{Xiaole Su, Hongwei Sun, Yusheng Wang\footnote{The
corresponding author (E-mail: wyusheng@bnu.edu.cn). \hfill{$\,$}}}

\vskip6mm

\noindent{\bf Abstract.} Let $M^n\ (n\geq3)$ be a complete
Riemannian manifold with $\sec_M\geq 1$, and let $M_i^{n_i}$
($i=1,2$) be two comlplete totally geodesic submanifolds in $M$. We
prove that if $n_1+n_2=n-2$ and if the distance
$|M_1M_2|\geq\frac{\pi}{2}$, then $M_i$ is isometric to
$\Bbb S^{n_i}/\Bbb Z_h$, ${\Bbb C\Bbb P}^{\frac {n_i}2}$
or ${\Bbb C\Bbb P}^{\frac {n_i}2}/\Bbb Z_2$ with the canonical
metric when $n_i>0$, and thus $M$ is isometric to $\Bbb
S^n/\Bbb Z_h$, ${\Bbb C\Bbb P}^{\frac n2}$ or ${\Bbb C\Bbb P}^{\frac
n2}/\Bbb Z_2$ except possibly when $n=3$ and $M_1$ (or $M_2$)
$\stackrel{\rm iso}{\cong}\Bbb S^{1}/\Bbb Z_h$ with $h\geq 2$ or $n=4$ and $M_1$
(or $M_2$) $\stackrel{\rm iso}{\cong}\Bbb{RP}^2$.

\vskip1mm


\vskip1mm


\vskip6mm

\setcounter{section}{-1}

\section{Introduction}

Let $M$ be a complete, simply connected Riemannian manifold with
$\sec_M\geq 1$.
$$\text{\it Under what conditions is $M$ isometric
to $\Bbb S^n$ or a projective space ${\Bbb K\Bbb P}^{n}$}?
\eqno{(0.1)}$$ Here, $\Bbb S^n$ is the unit sphere and ${\Bbb K\Bbb
P}^{n}$ is endowed with the canonical metric,
where $\Bbb K=\Bbb C, \Bbb H$ or $\Bbb Ca$ and $n\leq 2$ if
$\Bbb K=\Bbb Ca$, which satisfies
$$\sec_{{\Bbb S}^{n}}\equiv 1\text{ and } \diam(\Bbb S^n)=\pi, \text{
and } 1\leq\sec_{{\Bbb K\Bbb P}^{n}}\leq 4 \text{ and }\diam({\Bbb
K\Bbb P}^{n})=\frac\pi2.$$
This question draws lots of attention from geometrists.
Note that ``$\sec_M\geq 1$'' implies that the diameter $\diam(M)\leq\pi$.
Toponogov proved that
{\it if $\diam(M)=\pi$ (here that $M$ is simply connected is not needed),
then $M$ is isometric to $\Bbb S^n$} (Maximal Diameter Theorem).
And Berger proved that
{\it if $1\leq\sec_M\leq 4$, then either $\diam(M)>\frac\pi2$
and $M$ is homeomorphic to a sphere, or $\diam(M)=\frac\pi2$} and
$M$ is isometric to ${\Bbb S}^{n}(\frac12)$ or a ${\Bbb K\Bbb
P}^{n}$ (Minimal Diameter Theorem, [CE]). Afterwards, Grove-Shiohama
proved that {\it if $\diam(M)>\frac\pi2$ (here that $M$ is
simply connected is not needed), then $M$ is homeomorphic to a
sphere} ([GS]). Inspired by these, Gromoll-Grove, Wilhelm and Wilking proved step by step
that {\it if $\diam(M)=\frac\pi2$, then $M$ is either homeomorphic to a sphere,
or isometric to a ${\Bbb K\Bbb P}^{n}$} ($\frac\pi2$-Diameter
Rigidity Theorem, [GG1,2], [W], [Wi1]).

Note that the isometric classification in the $\frac\pi2$-Diameter
Rigidity Theorem is
on the premise that $M$ is not homeomorphic to a sphere. The present paper
aims to give ``purely isometric'' answers to question (0.1)
(as Toponogov's and Berger's above).

A basic fact is that ${\Bbb S}^{n}$ has a join structure, i.e.,
$$\Bbb S^n=\Bbb S^{n_1}*\Bbb S^{n_2},\eqno{(0.2)}$$
where $n_1+n_2=n-1$, and each $\Bbb S^{n_i}$ is totally geodesic in
$\Bbb S^n$, and the distance
$$|p_1p_2|=\frac\pi2\text{ for all } p_i\in\Bbb S^{n_i}.$$
Similarly, we can define a spherical join of two
Alexandrov spaces $X_i$ with curvature $\geq1$ (including
Riemannian manifolds with $\sec\geq 1$), $X_1*X_2$, which is also
an Alexandrov space with curvature $\geq1$ ([BGP]).  And, in $X_1*X_2$,
$X_i$ is convex\footnote{We say that a subset $A$ is convex in $X$
if, for any $x, y\in A$, there is a minimal geodesic (in $X$)
jointing $x$ with $y$ which belongs to $A$.}
with $\dim(X_1)+\dim(X_2)=\dim(X)-1$, and
$|p_1p_2|=\frac\pi2\text{ for all } p_i\in X_i$.

Inspired by this, Rong-Wang obtains the following rigidity theorem.

\vskip2mm

\noindent{\bf Theorem 0.1 ([RW])}{\bf.} {\it Let $X$ be a compact Alexandrov space
with curvature $\geq1$ and of dimension $n$, and let $X_i$ be its two compact convex subsets with empty
boundary and of dimension $n_i$.  If $|X_1X_2|\triangleq\min\{|p_1p_2||\ p_i\in X_i\}\geq\frac\pi2$,
then $n_1+n_2\leq n-1$, and equality implies that $X$
is isometric to a spherical join modulo a finite group.}

\vskip2mm

In Riemannian case we have the following corollary.

\vskip2mm

\noindent{\bf Corollary 0.2 ([RW])}{\bf.} {\it Let $M^n$ be a
complete Riemannian manifold with $\sec_M\geq 1$, and let
$M_i^{n_i}$ be its two complete totally geodesic submanifolds.
If $|M_1M_2|\geq\frac{\pi}{2}$, then $n_1+n_2\leq n-1$, and equality
implies that there is a finite group $\Gamma$ such that
$M_i$ is isometric to $\Bbb S^{n_i}/\Gamma$ and
$M$ is isometric to $\Bbb S^n/\Gamma$.}

\vskip2mm

Naturally, the next step is to consider the case where $n_1+n_2=n-2$.

\vskip2mm

\noindent{\bf Conjecture 0.3.} {\it
For $X$ and $X_i$ in Theorem 0.1, if
$|X_1X_2|\geq\frac\pi2$ and $n_1+n_2=n-2$, then either $X_i$ belong to a compact convex
subset of dimension $n-1$ in $X$, or $X$
is isometric to a spherical join modulo a $1$-dimensional Lie group (with finite components).}

\vskip2mm

\noindent{\bf Remark 0.4.} Note that this conjecture is trivial when $n_1$ or $n_2=0$.
The reason is that if $n_i=0$ and $X_i$ has an empty boundary, then it is our convention that
$X_i$ consists of two points with distance $\pi$. In this case,
$X$ is isometric to a spherical join, and isometric to a unit sphere
in Riemannian case (cf. [RW]).

\vskip2mm

In general case, Conjecture 0.3 will be much harder than Theorem 0.1.
In the present paper, the main result asserts that the conjecture is true
in Riemannian case.

First, let's focus on the complex projective space $\Bbb{CP}^{m}$
as an example satisfying the conjecture. In (0.2), we assume
that $n=2m+1$ and $n_i=2m_i+1$. We know that $S^1$ can act on $\Bbb
S^n$ freely and isometrically and preserve each $\Bbb
S^{n_i}$ such that
$$\Bbb S^n/S^1\stackrel{\rm iso}{\cong}\Bbb{CP}^{m}\text{ and } \Bbb S^{n_i}/S^1\stackrel{\rm
iso}{\cong}\Bbb{CP}^{m_i}. \eqno{(0.3)}$$
Note that $\Bbb{CP}^{m_i}$ is totally geodesic in $\Bbb{CP}^{m}$
with
$$|q_1q_2|=\frac\pi2\text{ for all } q_i\in\Bbb{CP}^{m_i},$$
and that
$\dim(\Bbb{CP}^{m_1})+\dim(\Bbb{CP}^{m_2})=\dim(\Bbb{CP}^{m})-2$.

Now, let's formulate our main result in this paper.

\vskip2mm

\noindent{\bf Main Theorem.}  {\it Let $M^n\ (n\geq 3)$ be a
complete Riemannian manifold with $\sec_M\geq 1$, and let
$M_i^{n_i}$ be its two complete totally geodesic submanifolds. If
$|M_1M_2|\geq\frac{\pi}{2}$ and $n_1+n_2=n-2$, then $M_i$ is
isometric to $\Bbb S^{n_i}/\Bbb Z_h$, ${\Bbb C\Bbb
P}^{\frac {n_i}2}$ or ${\Bbb C\Bbb P}^{\frac {n_i}2}/\Bbb Z_2$ when
$n_i>0$, and thus $M$ is isometric to $\Bbb S^n/\Bbb
Z_h$, ${\Bbb C\Bbb P}^{\frac n2}$ or ${\Bbb C\Bbb P}^{\frac n2}/\Bbb
Z_2$ except possibly when $n=3$ and $M_1$ (or $M_2$) is isometric to
$\Bbb S^{1}/\Bbb Z_h$ with $h\geq2$ or $n=4$ and $M_1$ (or $M_2$) is isometric to
$\Bbb{RP}^2$. }

\vskip2mm

Here, we make a convention that $M_i$ contains only one point if
$n_i=0$, and $\Bbb S^{1}$ is the circle with perimeter $2\pi$. And
one can refer to A.1 in Appendix for the construction of ${\Bbb
C\Bbb P}^{\frac n2}/\Bbb Z_2$ ($=\Bbb S^{n+1}/G$, where $G$ is a
$1$-dimensional Lie group with two components).

\vskip2mm

\noindent{\bf Remark 0.5.} Note that together with Remark 0.4, the Main Theorem implies
Conjecture 0.3 in Riemannian case.

\vskip2mm

\noindent{\bf Remark 0.6.} On the Main Theorem, we have the following
notes.

\noindent(0.6.1) If $\sec_{M_i}\not\equiv1$ ($n_i>1$) for $i=1$ or $2$
 or if there is an infinite number of minimal geodesics between some $p_1\in
M_1$ and $p_2\in M_2$, then $M\stackrel{\rm iso}{\cong}{\Bbb C\Bbb
P}^{\frac n2}$ or ${\Bbb C\Bbb P}^{\frac n2}/\Bbb Z_2$ (see
Proposition \ref{3.1} and Lemma \ref{3.14}), and that
$M\stackrel{\rm iso}{\cong}{\Bbb C\Bbb P}^{\frac n2}/\Bbb Z_2$
occurs only when $\frac {n_i}2$ and $\frac n2$ are all odd.

\noindent(0.6.2) If $n$ is odd, then $M\stackrel{\rm iso}{\cong}
\Bbb S^n/\Bbb Z_h$, and $h\geq3$ implies that $n_1$
or $n_2=0$; i.e., if in addition $n_1, n_2>0$, then $M\stackrel{\rm
iso}{\cong}\Bbb S^n$ or $\Bbb{RP}^n$ (see Lemma \ref{3.3} and
\ref{3.8}).

\noindent(0.6.3) Suppose that $M$ is simply connected. If $n\geq 5$, then $M\stackrel{\rm
iso}{\cong}\Bbb S^n$ or ${\Bbb C\Bbb P}^{\frac n2}$;
if $n=4$ (resp. $n=3$), either $M\stackrel{\rm iso}{\cong}\Bbb
S^4$ or ${\Bbb C\Bbb P}^2$ (resp. $\Bbb S^3$), or one of
$M_i\stackrel{\rm iso}{\cong}\Bbb{RP}^2$ (resp. $\Bbb S^{1}/\Bbb
Z_h$ with $h\geq 2$) (in this case, $M$ is homeomorphic to $\Bbb S^4$ (resp. $\Bbb
S^3$) by the $\frac\pi2$-Diameter Rigidity Theorem).

\vskip2mm

Based on (0.6.3), we have the following two questions.

\vskip2mm

\noindent{\bf Problem 0.7.} {\it Can ${\Bbb R\Bbb P}^2$
be embedded isometrically into $M^4$ as a totally geodesic submanifold,
where $M^4$ is a complete Riemannian manifold with $\sec\geq1$
and is homeomorphic to $\Bbb S^4$}?

\vskip2mm

\noindent{\bf Problem 0.8.} {\it In the Main Theorem, if $n\geq8$ and
$n_1+n_2\geq n-4$ (resp. $n\geq16$ and $n_1+n_2\geq n-8$),
is $M$ isometric to $\Bbb S^n$, ${\Bbb C\Bbb P}^{\frac n2}$ or
${\Bbb H\Bbb P}^{\frac n4}$ (resp. $\Bbb S^n$, ${\Bbb C\Bbb
P}^{\frac n2}$, ${\Bbb H\Bbb P}^{\frac n4}$ or ${\Bbb Ca\Bbb P}^2$)
when $M$ is simply connected}?

\vskip2mm

We can use the approach to the Main Theorem to discuss Problem 0.8,
but some essential difficulties will arise.

\vskip2mm

It seems that there is some overlap between the Main
Theorem and the $\frac\pi2$-Diameter Rigidity Theorem.
We will end this section by pointing out the main difference between them
by comparing the key points in their proofs.

\vskip2mm

\noindent{\bf Remark 0.9.}
(0.9.1) {\it The key point to the $\frac\pi2$-Diameter Rigidity Theorem}:

To the $\frac\pi2$-Diameter Rigidity Theorem,
an important fact is that $B'=B'''$ for any compact subset
$B\subset M$, where $B'=\{p\in M|\ |pb|=\frac\pi2\ \forall\ b\in B\}$
which is convex in $M$.
This is guaranteed by $\sec_M\geq1$ and $\diam(M)=\frac\pi2$
via Toponogov's Comparison Theorem (see Theorem 1.1 below).
In [GG1], $B'$ and $B''$ are called a pair of dual sets.
Then either both $B'$ and $B''$ are contractible, and in this case
$M$ is homeomorphic to a sphere; or both $B'$ and $B''$ are
totally geodesic submanifolds, and any $p\in M$ belongs to some
minimal geodesic $[q_1q_2]$ with $q_1\in B'$ and $q''\in B''$, and
then $M$ is isometric to a ${\Bbb K\Bbb P}^{n}$ (the proof involves
several big classification theorems ([GG2], [Wi1]), e.g. Bott-Samelson's
Theorem in [B]). In any case, the key fact that $B'$ and $B''$ are
dual to each other plays a crucial role.

\noindent (0.9.2) {\it The key point to the Main Theorem}:

In our Main Theorem, we in fact have that $|p_1p_2|=\frac\pi2$ for all $p_i\in M_i$
(see Corollary \ref{2.2} below). Hence, if $M_1$ and $M_2$ are dual to each other,
then we can use the approach to the $\frac\pi2$-Diameter Rigidity Theorem
to prove that $M$ is isometric to  ${\Bbb C\Bbb
P}^{\frac {n}2}$ or ${\Bbb C\Bbb P}^{\frac {n}2}/\Bbb Z_2$.
Indeed, if $n_1>0$ and $n_2>0$, then we can prove that either $M\stackrel{\rm
iso}{\cong}\Bbb S^n$ or $\Bbb{RP}^n$,
or $M_1$ and $M_2$ are dual to each other (see Proposition A.4
in Appendix). However, if one of $n_i=0$, we cannot see
that $M_1$ and $M_2$ are dual to each other. An important reason is that,
for any one of $M_i$ with $n_i>0$,
$$\left\{p\in M|\ |pM_i|\geq \frac\pi2\right\}=\left\{p\in M|\ |pp_i|=\frac\pi2\ \forall\ p_i\in M_i\right\}$$
(see Lemma \ref{2.1} below); but this may not be true if $n_i=0$
(i.e. $M_i$ is a single point). Therefore, the really challenging case to
the Main Theorem is where $M_1$ or $M_2$ is a single point.
Our proof for it, which also fits the case where $n_1>0$ and $n_2>0$,
is based on an easy observation that
$\lambda_{p_1p_2}=\lambda_{p_1'p_2'}$ for all $p_i,p_i'\in M_i$,
where $\lambda_{p_1p_2}$ denotes the number of all minimal
geodesics between $p_1$ and $p_2$ (see Corollary \ref{2.5} below).
If $\lambda_{p_1p_2}$ is finite (resp. infinite)
and $n_i>0$, then we can prove that $M_i$ is
isometric to $\Bbb S^{n_i}/\Bbb Z_h$ (resp. ${\Bbb C\Bbb
P}^{\frac {n_i}2}$ or ${\Bbb C\Bbb P}^{\frac {n_i}2}/\Bbb Z_2$).
In proving that  $M_i$ is isometric to ${\Bbb C\Bbb
P}^{\frac {n_i}2}$ or ${\Bbb C\Bbb P}^{\frac {n_i}2}/\Bbb Z_2$,
we do not use any big classification theorem involved in the proof of
the $\frac\pi2$-Diameter Rigidity Theorem, i.e.
our method is quite different from that in [GG1] and [Wi1].

%%%%%%%%%%%%%%%%%%%%%%%%%%%%%%%%%%%%%%% Section 1  %%%%%%%%%%%%%%%%%%%%%%%%%%%%%%%%%%%%%%%

\section{Toponogov's Comparison Theorem}

In this paper, we always let $[pq]$ denote a minimal geodesic
between $p$ and $q$ in a Riemannian manifold, and let $|pq|$
denote the distance between $p$ and $q$. Now, we give the main tool
of the paper---Toponogov's Comparison Theorem.

\begin{theorem}[{[P], [GM]}]\label{1.1}
Let $M$ be a complete Riemannian manifold with $\sec_M\geq\kappa$,
and let $\Bbb S^2_\kappa$ be the complete, simply
connected $2$-manifold of curvature $\kappa$.

\noindent {\rm(i)} To any $p\in M$ and $[qr]\subset M$,
we associate $\tilde p$ and $[\tilde q\tilde r]$
in $\Bbb S^2_\kappa$ with $|\tilde p\tilde q|=|pq|,|\tilde p\tilde
r|=|pr|$ and $|\tilde r\tilde q|=|rq|$. Then for any $s\in[qr]$ and
$\tilde s\in[\tilde q\tilde r]$ with $|qs|=|\tilde q\tilde s|$,
we have that $|ps|\geq|\tilde p\tilde s|$.

\vskip1mm

\noindent {\rm(ii)}  To any $[qp]$ and $[qr]$ in $M$,
we associate $[\tilde q\tilde p]$ and $[\tilde q\tilde r]$
in $\Bbb S^2_\kappa$ with $|\tilde q\tilde p|=|qp|$, $|\tilde q\tilde r|=|qr|$
and $\angle\tilde p\tilde q\tilde r=\angle pqr$. Then we have that
$|\tilde p\tilde r|\geq|pr|$.

\vskip1mm

\noindent {\rm(iii)} If the equality in {\rm (ii)} (resp. in {\rm (i)} for some $s$ in
the interior part of $[qr]$)
holds, then there exists a $[pr]$ (resp. $[pq]$ and $[pr]$) such that the triangle formed by $[qp]$, $[qr]$ and $[pr]$ bounds a
surface which is convex
and can be embedded isometrically into $\Bbb S^2_\kappa$.
\end{theorem}

%%%%%%%%%%%%%%%%%%%%%%%%%%%%%%%%%%%%%%% Section 2  %%%%%%%%%%%%%%%%%%%%%%%%%%%%%%%%%%%%%%%

\section{Preliminaries}

In this section, all $M_i\ (i=1,2)$ and $M$ are the manifolds in the
Main Theorem.

By (ii) of Theorem 1.1, one can prove the following interesting lemma.

\begin{lemma}[{[Ya]}]\label{2.1} Let $N$ be a complete Riemannian
manifold with $\sec_M\geq 1$, and let $L$ be a complete totally
geodesic submanifold in $N$ with $\dim(L)\geq 1$. Then $L^{\geq
\frac{\pi}{2}}=L^{=\frac{\pi}{2}}$.
\end{lemma}

In this lemma, $L^{\geq \frac{\pi}{2}}$ (resp. $L^{=\frac{\pi}{2}}$)
denotes the set $\{p\in N||px|\geq\frac{\pi}{2}\ \forall\ x\in L \}$
(resp. $\{p\in N||px|=\frac{\pi}{2}\ \forall\ x\in L \}$). (This
lemma has an Alexandrov version in [Ya], and one can refer to [SW]
for its detailed proof.) From Lemma \ref{2.1}, we can draw an
immediate corollary (which is fundamental and important to the
paper).

\begin{coro}\label{2.2} Under the conditions of the Main Theorem, we have that
$$|p_1p_2|=\frac{\pi}{2} \text{ for any } p_1\in M_1\text{ and } p_2\in M_2.\eqno{(2.1)}$$
\end{coro}

In this paper, we will let $\uparrow_{p}^{q}$ denote the unit tangent
vector at $p$ of a given (minimal geodesic) $[pq]$ (which is also
called the direction from $p$ to $q$ along $[pq]$); and let
$\Sigma_pM$ denote the set of all unit tangent vectors in $T_pM$. By
(2.1) and (ii) of Theorem 1.1, for any $[p_1p_2]$ and $[p_1p_2']$
with $p_1\in M_1$ and $p_2, p_2'\in M_2$, we have that
$$|\uparrow_{p_1}^{p_2}\uparrow_{p_1}^{p_2'}|\geq|p_2p_2'|.\eqno{(2.2)}$$

Now, we fix an arbitrary $[p_1p_2]$ with
$p_i\in M_i$. By Corollary \ref{2.2}, we conclude that
$$\uparrow_{p_2}^{p_1}\in (\Sigma_{p_2}M_2)^{=\frac{\pi}{2}}\subset \Sigma_{p_2}M.\eqno{(2.3)}$$
Similarly, we have that
$$\uparrow_{p_1}^{p_2}\in (\Sigma_{p_1}M_1)^{=\frac{\pi}{2}}\subset \Sigma_{p_1}M,\eqno{(2.4)}$$
where $(\Sigma_{p_1}M_1)^{=\frac{\pi}{2}}=\Sigma_{p_1}M$ if $n_1=0$.
By (iii) of Theorem 1.1, (2.3) and (2.1) imply the following easy fact.

\begin{lemma}\label{2.3} For any given $[p_1p_2]$ with $p_i\in M_i$
and $[p_2p_2']\subset M_2$, there exists a
$[p_1p_2']$ such that the triangle formed by $[p_1p_2]$,
$[p_2p_2']$ and $[p_1p_2']$ bounds a surface which is convex in $M$ and can be
embedded isometrically into the unit sphere $\Bbb S^2$.
\end{lemma}

For convenience, we call such a surface in Lemma \ref{2.3} a {\it
convex spherical surface}. Now, based on Lemma \ref{2.3}, we give
another important observation.

\begin{lemma}\label{2.4}
For any given $p_1\in M_1$ and $[p_2p_2']\subset M_2$,
there is a $1$-$1$ map
$$\iota:\{\text{minimal geodesics between $p_1$ and $p_2$}\}\to
\{\text{minimal geodesics between $p_1$ and $p_2'$}\}$$ such
that, for any $[p_1p_2]$, $\iota([p_1p_2])$ is the unique
minimal geodesic such that the triangle formed by $[p_1p_2]$, $\iota([p_1p_2])$
and $[p_2p_2']$ bounds a convex spherical surface.
\end{lemma}

Note that, in this lemma, if there is a sequence of minimal geodesics
$[p_1p_2]_j$ (between $p_1$ and $p_2$) with
$\lim\limits_{j\to\infty}[p_1p_2]_j\to[p_1p_2]$, then it is not hard to see
that $$\lim\limits_{j\to\infty}\iota([p_1p_2]_j)\to
\iota([p_1p_2]).\eqno{(2.5)}$$ And this lemma has an almost immediate
corollary.

\begin{coro}\label{2.5} Under the conditions of the Main Theorem, we have that
$$\lambda_{p_1p_2}=\lambda_{p_1'p_2'}\ \forall\ p_i,p_i'\in M_i,$$
where $\lambda_{p_1p_2}$ denotes the number of\ \{minimal
geodesics between $p_1$ and $p_2$\}.
\end{coro}

\noindent{\bf Proof of Lemma \ref{2.4}.} For any given
$[p_1p_2]$, by Lemma \ref{2.3}, there is a
$[p_1p_2']$ such that the triangle formed by $[p_1p_2]$, $[p_1p_2']$
and $[p_2p_2']$ bounds a convex spherical surface $D$. Note that $D$
determines a minimal geodesic $[\uparrow_{p_2}^{p_1}\uparrow_{p_2}^{p_2'}]$
of length $\frac\pi2$ in $\Sigma_{p_2}M$
(which is isometric to $\Bbb S^{n-1}$, so there is a unique minimal geodesic
between $\uparrow_{p_2}^{p_1}$ and $\uparrow_{p_2}^{p_2'}$).
Hence, for the given $[p_1p_2]$, such a $[p_1p_2']$ is
unique and vice versa, so the lemma follows. \hfill$\Box$

\vskip2mm

In the following, for any fixed $p_1\in M_1$, we will discuss the multi-valued
map  $$f_{p_1}:M_2\to (\Sigma_{p_1}M_1)^{=\frac{\pi}{2}} \text{ defined by }
p_2\mapsto\Uparrow_{p_1}^{p_2},$$
where $\Uparrow_{p_1}^{p_2}$ denotes the set of unit tangent vectors
at $p_1$ of all minimal geodesics between $p_1$ and $p_2$.
Note that $f_{p_1}$ is well defined because of (2.4). Obviously,
$(\Sigma_{p_1}M_1)^{=\frac{\pi}{2}}=\Bbb S^{n_2+1}$. For convenience, we let
$\Bbb S_{p_1}^{n_2+1}$ denote $(\Sigma_{p_1}M_1)^{=\frac{\pi}{2}}$.

We first note that, for any $[p_1p_2]$ with $p_i\in M_i$,
by Lemma \ref{2.4} we can define a map
$$f_{[p_1p_2]}:M_2\rightarrow \Bbb S_{p_1}^{n_2+1}$$
by $p_2\mapsto\uparrow_{p_1}^{p_2}$ and
$p_2'\mapsto\text{some}\uparrow_{p_1}^{p_2'}$ for any other $p_2'\in
M_2$ such that
$$|f_{[p_1p_2]}(p_2)f_{[p_1p_2]}(p_2')|=|p_2p_2'|,\eqno{(2.6)}$$
and that
$$f_{[p_1p_2]}([p_2p_2']) \text{\rm\ is a
minimal geodesic } [\uparrow_{p_1}^{p_2}\uparrow_{p_1}^{p_2'}] \text{\rm\ in }
\Bbb S_{p_1}^{n_2+1}\eqno{(2.7)}$$
if $[p_2p_2']$ is the unique minimal geodesic
between $p_2$ and $p_2'$. Then we can define a
``differential'' map  (cf. [RW])
$$\text{d}f_{[p_1p_2]}:\Sigma_{p_2}M_2\rightarrow
\Sigma_{\uparrow_{p_1}^{p_2}}\Bbb S_{p_1}^{n_2+1} \text{ by
}
\uparrow_{p_2}^{p_2'}\mapsto\uparrow_{\uparrow_{p_1}^{p_2}}^{\uparrow_{p_1}^{p_2'}}.$$
Note that $\Sigma_{p_2}M_2\stackrel{\text{iso}}{\cong}\Bbb
S^{n_2-1}$ and $\Sigma_{\uparrow_{p_1}^{p_2}}\Bbb
S_{p_1}^{n_2+1}\stackrel{\text{iso}}{\cong}\Bbb S^{n_2}$. About
$\text{\rm d}f_{[p_1p_2]}$, we have the following key observation.

\begin{lemma}\label{2.6} $\text{\rm d}f_{[p_1p_2]}$ is an isometrical embedding.
\end{lemma}

\noindent{\bf Proof.} By Lemma \ref{2.7} below, it suffices to show
that $\text{\rm d}f_{[p_1p_2]}$ is a distance nondecreasing map.
Note that (2.2) implies that
$$|f_{[p_1p_2]}(p_2')f_{[p_1p_2]}(p_2'')|\geq|p_2'p_2''|\eqno{(2.8)}$$
for all $p_2',p_2''\in M_2$. It is not hard to see that (2.8) together with (2.6)
and (2.7) implies that $\text{\rm d}f_{[p_1p_2]}$ is
a distance nondecreasing map. \hfill$\Box$

\begin{lemma}[{[SSW]}]\label{2.7}
Let $N$ be a complete Alexandrov space with curvature $\geq1$ (especially a complete
Riemannian manifold with $\sec_N\geq1$). If $f:\Bbb S^k\to N$ is a
distance nondecreasing map, then $f$ is an isometrical embedding.
\end{lemma}

Note that Lemma \ref{2.6} implies that there is an $\Bbb S^{n_2}$
passing $\uparrow_{p_1}^{p_2}$ in $\Bbb S_{p_1}^{n_2+1}$ such that
$$\Sigma_{\uparrow_{p_1}^{p_2}}\Bbb S^{n_2}=
\text{d}f_{[p_1p_2]}(\Sigma_{p_2}M_2).$$ For convenience, we let
$\Bbb S^{n_2}_{[p_1p_2]}$ denote this $\Bbb S^{n_2}$ (similarly, we
have the corresponding $\Bbb S^{n_1}_{[p_2p_1]}\
(\subset(\Sigma_{p_2}M_2)^{=\frac{\pi}{2}}=\Bbb S^{n_1+1})$ if $n_1>0$). By
the definition of $f_{[p_1p_2]}$ (together with Lemma \ref{2.4}), it
is not hard to see that
$$f_{[p_1p_2]}(M_2)\subseteq\Bbb S_{[p_1p_2]}^{n_2}.$$

\begin{remark}\label{2.8} {\rm By Lemma \ref{2.6},
$\text{\rm d}f_{[p_1p_2]}$ can be generalized naturally to an
isometry
$$\text{d}f_{[p_1p_2]}:T_{p_2}M_2\rightarrow
T_{\uparrow_{p_1}^{p_2}}\Bbb S_{[p_1p_2]}^{n_2}.$$ Then it is easy
to see that
$$f_{[p_1p_2]}|_{B_{M_2}(p_2,r_0)}=
\exp_{\uparrow_{p_1}^{p_2}}\circ\text{
d}f_{[p_1p_2]}\circ\left(\exp_{p_2}|_{B_{T_{p_2}M_2}(O,r_0)}\right)^{-1},$$ where
$r_0$ is the injective radius of $M_2$ (and $O$ is the original point of $T_{p_2}M_2$), and
$\exp_{\uparrow_{p_1}^{p_2}}$ and $\exp_{p_2}$ are the exponential
maps of $\Bbb S_{[p_1p_2]}^{n_2}$ and $M_2$ at
$\uparrow_{p_1}^{p_2}$ and $p_2$ respectively. (In the paper, we
denote by $B_{A}(p,r)$ the open $r$-ball in $A$ with the center $p$.) }
\end{remark}

Note that for any $[p_2p_2']\subset B_{M_2}(p_2,r_0)$,
$f_{[p_1p_2]}([p_2p_2'])$ is a minimal geodesic in $\Bbb
S_{[p_1p_2]}^{n_2}$ (see (2.7)), i.e., $f_{[p_1p_2]}([p_2p_2'])$ lies in a great
circle $\mathfrak{S}^1\subseteq \Bbb S^{n_2}_{[p_1p_2]}$. Let
$f_{[p_1p_2]}(p_2')$ be the direction of $[p_1p_2']$ (from $p_1$ to
$p_2'$). We can also consider the map $f_{[p_1p_2']}:M_2\to \Bbb
S_{p_1}^{n_2+1}$. Then it is not hard to conclude that:

\begin{lemma}\label{2.9} $\mathfrak{S}^1\subset f_{p_1}(M_2)$.
And $f_{p_1}^{-1}(\mathfrak{S}^1)$ is a closed geodesic
containing $[p_2p_2']$, and $f_{p_1}^{-1}|_{\mathfrak{S}^1}$ is a
locally isometrical map.
\end{lemma}

Lemma \ref{2.9} has the following almost immediate corollary.

\begin{coro}\label{2.10} {\rm (2.10.1)} For any $[p_1p_2]$
with $p_i\in M_i$, we have that $$\Bbb S^{n_2}_{[p_1p_2]}\subseteq
f_{p_1}(M_2).$$

\noindent{\rm (2.10.2)} Each minimal geodesic on $M_2$ lies in a
closed geodesic whose length divides $2\pi$.
\end{coro}

%%%%%%%%%%%%%%%%%%%%%%%%%%%%%%%%%%%%%%% Section 3  %%%%%%%%%%%%%%%%%%%%%%%%%%%%%%%%%%%%%%%

\section{Proof of The Main Theorem}

In this section, all $M_i\ (i=1,2)$ and $M$ are also the manifolds in the
Main Theorem, and we assume that $$n_2>0.$$ According to Corollary
\ref{2.5}, we can divide the whole proof into two parts: one is on
$\lambda_{p_1p_2}\equiv h<+\infty$, and the other is on
$\lambda_{p_1p_2}=+\infty$ for all $p_i\in M_i$. Hence, the Main
Theorem follows from Lemma \ref{3.3} and \ref{3.14} below (in the
proof of Lemma \ref{3.14}, Lemma \ref{3.10} plays the most important
role).

\vskip3mm

\noindent{\bf Part a.} On $\lambda_{p_1p_2}\equiv h<+\infty$ for all
$p_i\in M_i$.

\vskip3mm

We first give an observation on the condition
``$\lambda_{p_1p_2}\equiv h<+\infty$''.

\begin{prop}\label{3.1}
$\lambda_{p_1p_2}\equiv h<+\infty$ for all $p_i\in M_i$ if and only
if $n_2=1$ or $\sec_{M_2}\equiv1$.
\end{prop}

In its proof, we will use the classical Frankel's Theorem.

\begin{theorem}[{[Fr]}]\label{3.2} Let $N^n$ be a closed
positively curved manifold, and let $N_i^{n_i}\ (i=1,2)$  be
complete totally geodesic submanifolds in $N$. If $n_1+n_2\geq n$,
then $N_1\cap N_2\neq\emptyset$.
\end{theorem}

\noindent{\bf Proof of Proposition \ref{3.1}.}

If $\lambda_{p_1p_2}\equiv h<+\infty$ for all $p_i\in M_i$, the by
(3.4.1) below, for any $[p_1p_2]$ with $p_i\in M_i$, there is a
neighborhood $U\subset M_2$ of $p_2$ such that
$f_{[p_1p_2]}|_{U}:U\to \Bbb S^{n_2}_{[p_1p_2]}$ is an isometrical
embedding, which implies that $\sec_{M_2}\equiv1$ if $n_2>1$.

Next, it suffices to show that $\lambda_{p_1p_2}<+\infty$ (see
Corollary \ref{2.5}) by assuming that $n_2=1$ or
$\sec_{M_2}\equiv1$. We fix a $[p_1p_2]$ with $p_i\in M_i$, and
consider $\Bbb S^{n_2}_{[p_1p_2]}$. For any $\xi\in \Bbb
S^{n_2}_{[p_1p_2]}$, there is a $[p_1p_2']$ with $p_2'\in M_2$ such
that $\xi=\uparrow_{p_1}^{p_2'}$ (see (2.10.1)). {\bf Claim:}
$$\Bbb{S}_{[p_1p_2]}^{n_2}=\Bbb{S}_{[p_1p_2']}^{n_2}.$$ If $n_2=1$, then
$f_{p_1}^{-1}|_{\Bbb{S}_{[p_1p_2]}^1}:\Bbb{S}_{[p_1p_2]}^1\to M_2$
is a locally isometrical map (by Lemma \ref{2.9}), which implies the
claim. If $\sec_{M_2}\equiv1$, then by Remark \ref{2.8} we have that
$f_{[p_1p_2]}|_{B_{M_2}(p_2,\frac{r_0}{2})}:B_{M_2}(p_2,\frac{r_0}{2})\to
B_{\Bbb S^{n_2}_{[p_1p_2]}}(\uparrow_{p_1}^{p_2},\frac{r_0}{2})$ is
an isometry (where $r_0$ is the injective radius of $M_2$). It then
follows that for any $\eta\in B_{\Bbb
S^{n_2}_{[p_1p_2]}}(\uparrow_{p_1}^{p_2},\frac{r_0}{2})$ we have
that $\Bbb{S}_{[p_1p_2]}^{n_2}=\Bbb{S}_{[p_1p_2'']}^{n_2},$ where
$[p_1p_2'']$ with $p_2''\in M_2$ satisfies
$\uparrow_{p_1}^{p_2''}=\eta$. Then it is not hard to see that the
claim follows. On the other hand, by Theorem \ref{3.2} we have that
$$\Bbb S^{n_2}_{[p_1p_2]}\cap\Bbb S^{n_2}_{[p_1p_2]'}\neq\emptyset$$
for any other minimal geodesic $[p_1p_2]'$ between $p_1$ and $p_2$.
It then follows that
$$\Bbb S^{n_2}_{[p_1p_2]}=\Bbb S^{n_2}_{[p_1p_2]'}.$$ Hence,
$f_{p_1}(M_2)=\Bbb S^{n_2}_{[p_1p_2]}$, and
$f_{p_1}^{-1}:\Bbb S^{n_2}_{[p_1p_2]}\to M_2$ is a Riemannian covering map,
which implies that $\lambda_{p_1p_2}<+\infty$.
\hfill $\square$

\vskip2mm

Now we classify $M_i$ and $M$ under the condition
``$\lambda_{p_1p_2}\equiv h<+\infty$''.

\begin{lemma}\label{3.3}
If $\lambda_{p_1p_2}\equiv h<+\infty$ for all $p_i\in M_i$, then
{\rm (i)} when $n\geq5$, $M_i\stackrel{\rm iso}{\cong}\Bbb
S^{n_i}/\Bbb Z_h$ and $M\stackrel{\rm iso}{\cong}\Bbb S^{n}/\Bbb
Z_h$, and  $h\geq3$ implies that $n_1=0$; {\rm (ii)} when $n=4$,
$M_2\stackrel{\rm iso}{\cong}\Bbb{RP}^2$ and $M\stackrel{\rm
iso}{\cong}\Bbb{RP}^4$ (resp. $M_2\stackrel{\rm iso}{\cong}\Bbb
S^{2}$ and $M\stackrel{\rm iso}{\cong}\Bbb{S}^4$, or
$M_2\stackrel{\rm iso}{\cong}\Bbb{RP}^2$) if $M$ is not simply
connected (resp. $M$ is simply connected); {\rm (iii)} when $n=3$,
$M_2\stackrel{\rm iso}{\cong}\Bbb S^{1}$ and $M\stackrel{\rm
iso}{\cong}\Bbb{S}^3$, or $M_2\stackrel{\rm iso}{\cong}\Bbb
S^{1}/\Bbb Z_h$ with $h\geq2$.
\end{lemma}

In the proof of Lemma \ref{3.3}, we will use the following technical
lemmas. (For the convenience  of readers, we will give a brief proof
for (3.4.1) in Appendix.)

\begin{lemma}[{[RW]}]\label{3.4}
If $\lambda_{p_1p_2}\equiv h<+\infty$
for all $p_i\in M_i$, then for any $[p_1p_2]$

\noindent{\rm (3.4.1)} there is a neighborhood $U\subset M_2$ of
$p_2$ such that $f_{[p_1p_2]}|_{U}$ is an isometry;

\noindent{\rm (3.4.2)} there are neighborhoods $U_i\subset M_i$ of
$p_i$ such that $U_1*U_2$ \footnote{Refer to A.3 in Appendix for the
metric of $U_1*U_2$.} can be embedded isometrically into $M$ around
$[p_1p_2]$.
\end{lemma}

\begin{lemma}\label{3.5} Let $N^m$ be a complete
Riemannian manifold with $\sec_N\geq1$, and let $L^l$ be a complete
totally geodesic submanifold in $N$ with $l\geq\frac m2$. Assume
that $\sec_L\equiv1$.

\noindent{\rm (3.5.1)} If $l>\frac m2$, then we have that
$\sec_N\equiv1$; and if $m-l=2$ (resp. $m-l$ is odd) additionally,
then $\pi_1(N)=\pi_1(L)=\Bbb Z_k$ (resp. $\Bbb Z_2$ or $1$) for some
$k$;

\noindent{\rm (3.5.2)} If $l=\frac m2$, and if $L\stackrel{\text{\rm
iso}}{\cong}\Bbb{RP}^l$ and $N$ is not simply connected
additionally, then we have that $N\stackrel{\text{\rm
iso}}{\cong}\Bbb{RP}^m$ with the canonical metric.
\end{lemma}

In the proof of Lemma \ref{3.5}, we will use the following
connectedness theorem.

\begin{theorem}[{[Wi2]}] \label{3.6} Let $N^m$ be a closed
positively curved manifold, and let $L^l$ be a complete totally
geodesic submanifold in $N$ with $l\geq\frac m2$. Then
$L\hookrightarrow N$ is $(2l-m+1)$-connected.
\end{theorem}

\noindent{\bf Proof of Lemma \ref{3.5}.}

Let $\pi:\tilde N\to N$ be the Riemannian covering map, and let
$\tilde L=\pi^{-1}(L)$ which is complete and totally geodesic in
$\tilde N$. By Theorem \ref{3.2}, $\tilde L$ is connected because
$l\geq\frac m2$.

(3.5.1) Since $l>\frac m2$, both $L\hookrightarrow N$ and $\tilde
L\hookrightarrow\tilde N$ are at least 2-connected (Theorem
\ref{3.6}). This implies that $\pi_1(N)=\pi_1(L)$ and $\tilde L$ is
simply connected. It then follows that $\tilde
L\stackrel{\text{iso}}{\cong}\Bbb{S}^l$ (note that $\sec_L\equiv1$).
By the Maximal Diameter Theorem, it has to hold that $\tilde
N\stackrel{\text{iso}}{\cong}\Bbb{S}^m$, i.e. $\sec_N\equiv1$. And
it is easy to see that $\pi_1(N)=\pi_1(L)=\Bbb Z_2$ or $1$ if $m-l$
is odd. Now we assume that $m-l=2$. Note that there is a great
circle $\Bbb S^1$ such that $\tilde N=\tilde L*\Bbb{S}^1$ (i.e.
$\Bbb{S}^m=\Bbb S^l*\Bbb S^1$). On the other hand, $\pi_1(N)\
(=\pi_1(L))$ acts on $\tilde N$ freely by isometries. Moreover,
$\pi_1(N)$ preserves $\tilde L$, and thus it also preserves the
$\Bbb S^1$. It follows that $\pi_1(N)=\pi_1(L)=\Bbb Z_k$ for some
$k$.

(3.5.2)  Since $m=2l$, we know that $\pi_1(N)=\Bbb Z_2$. On the
other hand, since $\tilde L=\pi^{-1}(L)$ (which is connected) and
$L\stackrel{\text{\rm iso}}{\cong}\Bbb{RP}^l$, it has to hold that
$\tilde L\stackrel{\text{iso}}{\cong}\Bbb{S}^l$. Similarly, by the
Maximal Diameter Theorem we get that $\tilde
N\stackrel{\text{iso}}{\cong}\Bbb{S}^m$, and so
$N\stackrel{\text{iso}}{\cong}\Bbb{RP}^m$. \hfill $\square$

\vskip2mm

Now we give the proof of Lemma \ref{3.3}.

\vskip2mm

\noindent{\bf Proof of Lemma \ref{3.3}.}

From the proof of Proposition \ref{3.1}, for a fixed $[p_1p_2]$ with
$p_i\in M_i$, $f_{p_1}^{-1}:\Bbb S^{n_2}_{[p_1p_2]}\to M_2$ is a
Riemannian covering map, and \#$(\pi_1(M_2))=h$ when $n_2\geq 2$.
Next, we will divide the proof into
the following two cases.

Case 1: $n_1=0$. In this case, $n_2=n-2$. If $n\geq5$, then by
(3.5.1) we have that $M_2\stackrel{\rm iso}{\cong}\Bbb S^{n_2}/\Bbb
Z_h$, and $M\stackrel{\rm iso}{\cong}\Bbb S^{n}/\Bbb Z_h$ (note that
$n_2>\frac n2$ and \#$(\pi_1(M_2))=h$). If
$n=4$, then $M_2\stackrel{\rm iso}{\cong}\Bbb S^2$ or $\Bbb{RP}^2$
(note that $n_2=2$), and thus respectively,
$M\stackrel{\text{iso}}{\cong}\Bbb{S}^4$ by the Maximum Diameter
Theorem or $M\stackrel{\text{iso}}{\cong}\Bbb{RP}^4$ by (3.5.2) if
$M$ is not simply connected. If $n=3$, then $M_2\stackrel{\rm
iso}{\cong}\Bbb S^1/\Bbb Z_h$ (note that $n_2=1$), and thus
$M\stackrel{\text{iso}}{\cong}\Bbb{S}^3$ by the Maximum Diameter
Theorem if $M_2\stackrel{\rm iso}{\cong}\Bbb S^1$.

Case 2: $n_1>0$. In this case, for any $p_i\in M_i$, we have proved
that $f_{p_1}(M_2)=\Bbb S^{n_2}$ and $f_{p_2}(M_1)=\Bbb S^{n_1}$,
and both $f_{p_1}^{-1}:\Bbb S^{n_2}\to M_2$ and $f_{p_2}^{-1}:\Bbb
S^{n_1}\to M_1$ are Riemannian covering maps. Together with (3.4.2),
this implies that the set
$$N\triangleq\{p\in M|\ p\text{ belongs to some $[p_1p_2]$
with $p_i\in M_i$}\}$$ is a complete totally geodesic
$(n-1)$-dimensional submanifold in $M$ (Hint: It follows from Lemma
\ref{2.3} that, for any $p_i\in M_i$, say $p_1$,
$\Sigma_{p_1}N=f_{p_1}(M_2)*\Sigma_{p_1}M_1=\Bbb S^{n_2}*\Bbb
S^{n_1-1}=\Bbb S^{n-2}$). Hence, by Corollary 0.2, we know that
$\sec_N\equiv1$ (note that $M_1$ and $M_2$ are totally geodesic in
$N$). It then follows from (3.5.1) that $M$ is isometric to $\Bbb
S^n$ or $\Bbb R\Bbb P^n$ (and so $M_i$ is isometric to $\Bbb
S^{n_i}$ or $\Bbb R\Bbb P^{n_i}$ respectively). \hfill$\Box$

\begin{remark}\label{3.7}{\rm Why cannot we prove that
$M\stackrel{\text{iso}}{\cong}\Bbb{RP}^4$ (resp. $M\stackrel{\rm
iso}{\cong}\Bbb S^3/\Bbb Z_h$ with $h\geq 2$)
when $n=4$ and $M_2\stackrel{\text{iso}}{\cong}\Bbb{RP}^2$
(resp. $n=3$ and $M_2\stackrel{\text{iso}}{\cong}\Bbb S^1/\Bbb Z_h$)
by a similar argument to the above proof for $n_1>0$?
Note that in such two cases, $n_1=0$, i.e. $M_1=\{p_1\}$.
Due to the similarity, we only give an explanation for the case where $n=4$.
Note that $\lambda_{p_1p_2}\equiv 2$ in this case, i.e., there are only two
minimal geodesics $[p_1p_2]_j$ ($j=1,2$) between $p_1$ and any $p_2\in M_2$.
Since $f_{p_1}^{-1}:\Bbb S^{2}_{[p_1p_2]_j}\to M_2$ is a Riemannian covering map,
$[p_1p_2]_1$ and $[p_1p_2]_2$ form an angle equal to $\pi$ at
$p_1$. However, we cannot judge
whether they form an angle equal to $\pi$
at $p_2$ or not, so that we cannot judge whether $N\triangleq\{p\in M|\ p\text{ belongs to some $[p_1p_2]$
with $p_2\in M_2$}\}$ is totally geodesic
in $M$ or not. }
\end{remark}

\vskip1mm

\noindent{\bf Part b.} On $\lambda_{p_1p_2}=+\infty$ for all $p_i\in M_i$.

\vskip3mm

\begin{lemma}\label{3.8}
If $\lambda_{p_1p_2}=+\infty$ for all $p_i\in M_i$, then both $n_1$
and $n_2$ are even.
\end{lemma}

\noindent{\bf Proof.} By Proposition \ref{3.1}, it suffices to
derive a contradiction by assuming that $n_2=2m+1$ with $m>0$. We
still fix an arbitrary $[p_1p_2]$ with $p_i\in M_i$ at first. And,
in this proof, we always let $\tilde q$ denote $f_{[p_1p_2]}(q)$ for
any $q\in B_{M_2}(p_2,\frac{r_0}2)$, where $r_0$ is the injective
radius of $M_2$.

{\bf Claim 1}: {\it There is an $\Bbb S^{m+1}\subset\Bbb
S^{2m+1}_{[p_1p_2]}$ such that $f_{p_1}^{-1}|_U$ is an isometry for some
convex domain $U$ in the $\Bbb S^{m+1}$}. We will find such an $\Bbb
S^{m+1}$ through the following steps.

\vskip1mm

Step 1. We select an arbitrary $\tilde p_2^1\in B_{\Bbb
S^{2m+1}_{[p_1p_2]}}(\tilde
p_2,\frac{r_0}{2})\setminus\{\tilde p_2\}$. By the definition of
$f_{[p_1p_2]}$, there is a $[p_1p_2^1]$ with
$\uparrow_{p_1}^{p_2^1}=\tilde p_2^1$ such that
$f_{p_1}^{-1}|_{[\tilde p_2\tilde p_2^1]}: [\tilde p_2\tilde
p_2^1]\to[p_2p_2^1]$ is an isometry (see (2.7)). For convenience, we
denote by $\Bbb S^1_{\bullet}$ the great circle including $[\tilde
p_2\tilde p_2^1]$ in  $S^{2m+1}_{[p_1p_2]}$.

We also consider $\Bbb S^{2m+1}_{[p_1p_2^1]}$, and
observe that
$$[\tilde p_2\tilde p_2^1]\subset\Bbb S^{2m+1}_{[p_1p_2]}\cap\Bbb
S^{2m+1}_{[p_1p_2^1]}$$ and
$$\Bbb S^{2m+1}_{[p_1p_2]}\cap\Bbb
S^{2m+1}_{[p_1p_2^1]}=\Bbb S^{k_1}\ (\text{denoted by } \Bbb
S^{k_1}_{[p_1p_2]}) \text{ with } k_1\geq 2m\eqno{(3.1)}$$ (note
that $\Bbb S^{2m+1}_{[p_1p_2]}, \Bbb
S^{2m+1}_{[p_1p_2^1]}\subset\Bbb S_{p_1}^{2m+2})$. Note that
$$f_{p_1}^{-1}(B_{\Bbb S^{k_1}_{[p_1p_2]}}(\tilde p_2,\frac{r_0}{2}))
\subset B_{M_2}(p_2,\frac{r_0}{2})\cap B_{M_2}(p_2^1,r_0),$$ and
thus, for any $\tilde p_2'\in B_{\Bbb S^{k_1}_{[p_1p_2]}}(\tilde
p_2,\frac{r_0}{2})$, we have that
$$|\tilde p_2\tilde p_2^1|=|p_2p_2^1|, |\tilde p_2\tilde p_2'|=|p_2p_2'|,
|\tilde p_2^1\tilde p_2'|=|p_2^1p_2'|.$$ Moreover, note that $\angle
p_2'p_2p_2^1=\angle \tilde p_2'\tilde p_2\tilde p_2^1$ (Lemma
\ref{2.6} and Remark \ref{2.8}). It then follows from (iii) of
Theorem 1.1 that
$$\text{\it the triangle $\triangle p_2p_2^1p_2'$ bounds a convex spherical surface in
$M_2$},\eqno{(3.2)}$$
where the triangle $\triangle p_2p_2^1p_2'$ is formed by $[p_2p_2^1]$,
$[p_2p_2']$ and $[p_2^1p_2']$ (note that there is a unique
minimal geodesic between any two points in
$B_{M_2}(p_2,\frac{r_0}{2})$).

\vskip1mm

Step 2. We select $\tilde p_2^2\in B_{\Bbb
S^{k_1}_{[p_1p_2]}}(\tilde p_2,\frac{r_0}{2})\setminus\Bbb
S^1_{\bullet}$, and let $[p_1p_2^2]$ be the minimal geodesic such that $\tilde
p_2^2=\uparrow_{p_1}^{p_2^2}$. And we let $\Bbb S^2_{\bullet}$ be
the unit sphere $\Bbb S^2\subset \Bbb S^{k_1}_{[p_1p_2]}$ including
$\tilde p_2, \tilde p_2^1, \tilde p_2^2$,  and let $D$ be the convex
domain in $\Bbb S^2_{\bullet}$ bounded by $\triangle \tilde
p_2\tilde p_2^1\tilde p_2^2$. By (3.2), it is easy to see that
$f_{p_1}^{-1}|_D$ is an isometry.

Similarly, we consider
$$\Bbb S^{k_1}_{[p_1p_2]}\cap\Bbb
S^{2m+1}_{[p_1p_2^2]}=\Bbb S^{k_2}\ (\text{denoted by } \Bbb
S^{k_2}_{[p_1p_2]}) \text{ with } k_2\geq 2m-1;$$ and, for any $\tilde
p_2'\in B_{\Bbb S^{k_2}_{[p_1p_2]}}(\tilde p_2,\frac{r_0}{2})$,
by (iii) of Theorem 1.1, we can derive that
$$\{p_2',p_2, p_2^1, p_2^2\} \text{\ as vertices
determines a convex spherical tetrahedron in } M_2.\eqno{(3.3)}$$

$\cdots$

Step $m+1$. We select $\tilde p_2^{m+1}\in B_{\Bbb
S^{k_{m}}_{[p_1p_2]}}(\tilde p_2,\frac{r_0}{2})\setminus\Bbb
S^m_{\bullet}$. Let $\Bbb S^{m+1}_{\bullet}$ be the unit  sphere
$\Bbb S^{m+1}\subset \Bbb S^{k_m}_{[p_1p_2]}$ including $\tilde p_2,
\tilde p_2^1, \cdots, \tilde p_2^{m+1}$ which as vertices determines
a convex domain $U$ in $\Bbb S^{m+1}_{\bullet}$. Similarly, by the
corresponding property similar to (3.2) and (3.3) in the $m$-th
step, we have that $f_{p_1}^{-1}|_U$ is an isometry. That is, $\Bbb
S^{m+1}_{\bullet}$ is just the wanted sphere in Claim 1.

\vskip1mm

In fact, Claim 1 has the following strengthened version.

\noindent{\bf Claim 2}: {\it For any $\tilde p_2'\in \Bbb
S^{m+1}_{\bullet}$, $f_{p_1}^{-1}|_{B_{\Bbb
S^{m+1}_{\bullet}}(\tilde p_2',\frac{r_0}{2})}$ is an isometry}. We
first select a $\tilde p_{2,0}'$ in the interior part of $U\subset
\Bbb S^{m+1}_{\bullet}$. Let $p_{2,0}'=f_{p_1}^{-1}(\tilde
p_{2,0}')$ and $[p_1p_{2,0}']$ be the minimal geodesic such that
$\tilde p_{2,0}'=\uparrow_{p_1}^{p_{2,0}'}$. Since $f_{p_1}^{-1}|_U$
is an isometry, it is easy to see that $\Bbb
S^{m+1}_{\bullet}\subset \Bbb S^{2m+1}_{[p_1p_{2,0}']}\cap \Bbb
S^{2m+1}_{[p_1p_2'']}$, where $[p_1p_2'']$ is the minimal geodesic
such that $\uparrow_{p_1}^{p_2''}$ belongs to $U$ and is close to
$\tilde p_{2,0}'$. By the arguments to get (3.2) we can conclude
that $f_{p_1}^{-1}|_{B_{\Bbb S^{m+1}_{\bullet}}(\tilde
p_{2,0}',\frac{r_0}{2})}$ is an isometry. Then by replacing $U$ with
$B_{\Bbb S^{m+1}_{\bullet}}(\tilde p_{2,0}',\frac{r_0}{2})$, it is
not hard to see that $f_{p_1}^{-1}|_{B_{\Bbb
S^{m+1}_{\bullet}}(\tilde p_2',\frac{r_0}{2})}$ is an isometry for
any $\tilde p_2'\in \Bbb S^{m+1}_{\bullet}$.

\vskip1mm

Inspired by the proof of Claim 2, we have the following observation.

\noindent{\bf Claim 3}: {\it For any small $\varepsilon>0$, there is
another minimal geodesic $[p_1p_2]'$ between $p_1$ and $p_2$ such
that $|\uparrow_{p_1}^{p_2}(\uparrow_{p_1}^{p_2})'|<\varepsilon$}.
In fact, if this is not true, then based on (2.5) we can use a
similar proof of (3.4.1) (ref. A.2 in Appendix) to prove that there
is a neighborhood $V\subset M_2$ of $p_2$ such that
$f_{[p_1p_2]}|_{V}$ is an isometry, so is
$f_{p_1}^{-1}|_{f_{[p_1p_2]}(V)}$. Note that $f_{[p_1p_2]}(V)$ is an
open subset in $\Bbb S^{2m+1}_{[p_1p_2]}$, so the proof of Claim 2
implies that $f_{p_1}^{-1}|_{B_{\Bbb S^{2m+1}_{[p_1p_2]}}(\tilde
p_2',\frac{r_0}{2})}$ is an isometry for any $\tilde p_2'\in \Bbb
S^{2m+1}_{[p_1p_2]}$. This implies that $f_{p_1}^{-1}|_{\Bbb
S^{2m+1}_{[p_1p_2]}}$ is a Riemannian covering map, which
contradicts Proposition \ref{3.1}.

\vskip1mm

Now, we will complete the whole proof of the lemma based on Claims
1-3. We still consider the above $\Bbb S^{m+1}_{\bullet}
\subset\Bbb S^{2m+1}_{[p_1p_2]}$ (with $\uparrow_{p_1}^{p_2}\in \Bbb S^{m+1}_{\bullet}$).
Let $[p_1p_2]'$ be another minimal
geodesic between $p_1$ and $p_2$ with $(\uparrow_{p_1}^{p_2})'$
being sufficiently close to $\uparrow_{p_1}^{p_2}$. We consider the
natural isometry (ref. Remark 2.8)
$$h\triangleq \exp_{(\uparrow_{p_1}^{p_2})'}\circ\text{ d}f_{[p_1p_2]'}
\circ\text{d}f_{[p_1p_2]}^{-1}\circ\exp_{\uparrow_{p_1}^{p_2}}^{-1}:\Bbb
S^{2m+1}_{[p_1p_2]}\to\Bbb S^{2m+1}_{[p_1p_2]'}.$$ Duo to (2.5) and
that $(\uparrow_{p_1}^{p_2})'$ is sufficiently close to
$\uparrow_{p_1}^{p_2}$, it is easy to see that $$|h(\tilde
p_2')\tilde p_2'|\ll\frac{r_0}{2}\ (\text{in $\Bbb S^{2m+2}_{p_1}$)
for any $\tilde p_2'\in \Bbb S^{2m+1}_{[p_1p_2]}$}.$$ On the other
hand, it is not hard to see that the unit sphere $h(\Bbb
S^{m+1}_{\bullet})$ (containing $(\uparrow_{p_1}^{p_2})'$) satisfies
that $f_{p_1}^{-1}|_{B_{h(\Bbb S^{m+1}_{\bullet})}(\tilde
p_2'',\frac{r_0}{2})}$ is an isometry for any $\tilde p_2''\in
h(\Bbb S_{\bullet}^{m+1})$. By Theorem \ref{3.2}, we have that $\Bbb
S^{m+1}_{\bullet}\cap h(\Bbb S^{m+1}_{\bullet})\neq\emptyset$ in
$\Bbb S_{p_1}^{2m+2}$. Select a $\tilde q$ in $\Bbb
S^{m+1}_{\bullet}\cap h(\Bbb S^{m+1}_{\bullet})$, and let
$q=f_{p_1}^{-1}(\tilde q)$. Since $f_{p_1}^{-1}|_{B_{h(\Bbb
S^{m+1}_{\bullet})}(\tilde q,\frac{r_0}{2})}$ is an isometry and
$|h(\tilde q)\tilde q|<\frac{r_0}{2}$, it has to hold that $h(\tilde
q)=\tilde q$ (note that $\tilde q, h(\tilde q)\in f_{p_1}(q)$),
which implies that $f_{p_1}^{-1}([\uparrow_{p_1}^{p_2}\tilde q])$
and $f_{p_1}^{-1}([(\uparrow_{p_1}^{p_2})'\tilde q])$ are the same
geodesic between $p_2$ and $q$ in $M_2$. This contradicts Lemma
\ref{2.6} once $f_{[p_1q]}$ is considered, where $[p_1q]$ is the
minimal geodesic such that $\uparrow_{p_1}^{q}=\tilde q$.
\hfill$\Box$

\begin{remark}\label{3.9}{\rm
If $n_2=2m+2$ in the above proof, then $\Bbb S^{m+1}_{\bullet}\cap
h(\Bbb S^{m+1}_{\bullet})$ can be empty in $\Bbb S_{p_1}^{2m+3}$
(and we can not find an $\Bbb S^{m+2}\subset\Bbb S_{[p_1p_2]}^{2m+2}$
such that $f_{p_1}^{-1}|_{\Bbb S^{m+2}}$ is a local isometry). In
fact, if $M\stackrel{\rm iso}{\cong}\Bbb {CP}^{\frac n2}$ and
$M_2\stackrel{\rm iso}{\cong}\Bbb {CP}^{m+1}$, then
$f_{p_1}^{-1}(\Bbb S^{m+1}_\bullet)$ is a complete totally geodesic
submanifold (in $M_2$) which is isometric to $\Bbb{RP}^{m+1}$ with the
canonical metric. }
\end{remark}

By Lemma \ref{3.8}, we can assume that $n_2=2m>0$.

\begin{lemma}\label{3.10}
If $\lambda_{p_1p_2}=+\infty$ for all $p_i\in M_i$, then $M_2$ is
isometric to ${\Bbb C\Bbb P}^{m}$ or ${\Bbb C\Bbb P}^{m}/\Bbb Z_2$
with the canonical metric, and that $M_2\stackrel{\rm iso}{\cong}\Bbb
{CP}^{m}/\Bbb Z_2$ occurs only when $m$ is odd.
\end{lemma}

In order to prove Lemma \ref{3.10}, we first give a key observation.

\begin{lemma}\label{3.11} For any $\epsilon>0$, there is a $\delta>0$
such that if
$|(\uparrow_{p_1}^{p_2})_1(\uparrow_{p_1}^{p_2})_2|<\delta$ for
arbitrary two minimal geodesics $[p_1p_2]_1$ and $[p_1p_2]_2$
between $p_1\in M_1$ and $p_2\in M_2$, then
$$\left||\uparrow_{\tilde p_2^1}^{\tilde p_2^2}
\xi|-\frac{\pi}{2}\right| <\epsilon$$ for any $\xi\in\Sigma_{\tilde
p_2^1}\Bbb S_{[p_1p_2]_1}^{2m}$, where $\tilde p_2^j$ denotes
$(\uparrow_{p_1}^{p_2})_j$.

\end{lemma}

\noindent{\bf Proof.} On $\Bbb S_{p_1}^{2m+1}$, for any $p_2'\in
B_{M_2}(p_2,r_0)$ (where $r_0$ is the injective radius of $M_2$),
$$f_{[p_1p_2]_1}(p_2')\in \Bbb S_{[p_1p_2]_1}^{2m} \text{ and }
|\tilde p_2^1f_{[p_1p_2]_1}(p_2')|=|p_2p_2'|,$$ and by (ii) of
Theorem \ref{1.1}
$$|\tilde p_2^2f_{[p_1p_2]_1}(p_2')|\geq|p_2p_2'|.$$
It then is easy to see that the lemma follows from the first
variation formula. \hfill$\Box$

\vskip2mm

\noindent{\bf Proof of Lemma \ref{3.10}.}

We still fix an arbitrary
$[p_1p_2]$ with $p_i\in M_i$ at first, and consider $f_{p_1},
f_{[p_1p_2]}, \Bbb S_{p_1}^{2m+1}, \Bbb S_{[p_1p_2]}^{2m}$ and so
on. By (2.7), for any $[p_2p_2']\subset B_{M_2}(p_2,\frac{r_0}{2})$,
there is a $[p_1p_2']$ such that
$f_{[p_1p_2]}|_{[p_2p_2']}:[p_2p_2']\to[\uparrow_{p_1}^{p_2}\uparrow_{p_1}^{p_2'}]$
is an isometry. Similar to the proof of Lemma \ref{3.8}, we have
that
$$[\uparrow_{p_1}^{p_2}\uparrow_{p_1}^{p_2'}]\subset\Bbb S^{2m}_{[p_1p_2]}\cap\Bbb
S^{2m}_{[p_1p_2']},$$ and we will consider (similar to (3.1))
$$\Bbb S^{k_1}_{[p_1p_2],[p_2p_2']}\triangleq\Bbb S^{2m}_{[p_1p_2]}\cap\Bbb
S^{2m}_{[p_1p_2']} \text{ with } k_1\geq 2m-1.\eqno{(3.4)}$$ {\bf
Claim 1}: {\it In fact, we have that $k_1=2m-1$}. If $k_1=2m$, then
similar to Claims 1 and 2 in the proof of Lemma \ref{3.8} we can
find an $\Bbb S^{m+1}\subset\Bbb S_{p_1}^{2m+1}$ such that
$f_{p_1}^{-1}|_{\Bbb S^{m+1}}$ is a local isometry, and we can
obtain a contradiction.

\vskip1mm

For convenience, we let $\gamma(t)|_{t\in[0,|p_2p_2'|]}$ denote
the $[p_2p_2']$ with $\gamma(0)=p_2$ ($t$ is the arc-length parameter), and let
$\tilde\gamma(t)$ denote $f_{[p_1p_2]}(\gamma(t))$. (In this proof,
we also let $\tilde q$ denote $f_{[p_1p_2]}(q)$ for any $q\in
B_{M_2}(p_2,\frac{r_0}2)$). Let $[p_1\gamma(t)]$ be the minimal geodesic
such that $\uparrow_{p_1}^{\gamma(t)}=\tilde\gamma(t)$. We claim
that
$$\Bbb S^{2m}_{[p_1\gamma(t)]}\cap\Bbb
S^{2m}_{[p_1\gamma(t')]}=\Bbb S^{2m-1}_{[p_1p_2],[p_2p_2']} \text{
for all } t\neq t'.\eqno{(3.5)}$$ We need only to verify it for
$m\geq2$. Note that for any $\Bbb S^2\subset\Bbb
S^{2m-1}_{[p_1p_2],[p_2p_2']}$ containing $[\tilde p_2\tilde p_2']$,
$f_{p_1}^{-1}|_{B_{\Bbb S^2}(\tilde p_2,\frac{r_0}2)}$ is an
isometry (which is similar to ``$f_{p_1}^{-1}|_{B_{\Bbb
S^{m+1}_{\bullet}}(\tilde p_2,\frac{r_0}2)}$ is an isometry'' in Claim 2
in the proof of Lemma \ref{3.8}, and implies that
$f_{p_1}^{-1}(B_{\Bbb S^2}(\tilde p_2,\frac{r_0}2))$ is totally geodesic
in $M_2$). This implies that $\Bbb S^2\subset \Bbb S^{2m}_{[p_1\gamma(t)]}\cap\Bbb
S^{2m}_{[p_1\gamma(t')]}$, and so (3.5) follows from Claim 1. Moreover, a parallel
vector field $X(t)$ along $\gamma(t)$ on $f_{p_1}^{-1}(B_{\Bbb S^2}(\tilde p_2,\frac{r_0}2))$ is
also parallel on $M_2$, and we can naturally define
$\text{d}f_{[p_1p_2]}(X(t))$, denoted by $\tilde X(t)$, which is
parallel along $\tilde\gamma(t)$ on $\Bbb S^2$ and satisfies
$|\tilde X(t)|=|X(t)|$. Note that we can select such $2m-2$ parallel
and orthogonal (unit) vector fields $X_1(t),\cdots,X_{2m-2}(t)$
along $\gamma(t)$ which are all perpendicular to $\gamma'(t)$, the tangent
vector of $\gamma(t)$. And
by Lemma \ref{2.6}, $\tilde X_1(t),\cdots,\tilde X_{2m-2}(t)$
are also orthogonal and perpendicular to $\tilde\gamma'(t)$.

Now we select a parallel unit vector field $X_{2m-1}(t)$ along
$\gamma(t)$ (on $M_2$) which is perpendicular to $X_1(t),\cdots,X_{2m-2}(t)$
and $\gamma'(t)$.

\noindent{\bf Claim 2}: {\it We can define $\text{\rm
d}f_{[p_1p_2]}(X_{2m-1}(t))$, denoted by $\tilde X_{2m-1}(t)$, which
is smooth with respect to $t$ and perpendicular to $\tilde
X_1(t),\cdots,$ $\tilde X_{2m-2}(t)$ and $\tilde\gamma'(t)$, and
satisfies $|\tilde X_{2m-1}(t)|\geq1.$} Let
$\beta_t(s)|_{s\in[0,\epsilon)}\subset B_{M_2}(p_2,\frac{r_0}2)$ be
a geodesic such that $\beta_t(0)=\gamma(t)$ and
$\beta_t'(0)=X_{2m-1}(t)$. Due to Remark \ref{2.8}, each
$\tilde\beta_t(s)|_{s\in[0,\epsilon)}$ is a smooth curve with
respect to $s$ (but it will not be a geodesic when $t>0$). It then
follows that we can define $$\text{\rm
d}f_{[p_1p_2]}(X_{2m-1}(t))\triangleq\tilde\beta_t'(s)|_{s=0}\
(\text{denoted by }\tilde X_{2m-1}(t)),$$ which is smooth with
respect to $t$ (this is also due to Remark \ref{2.8}). On the other
hand, since $|\tilde p_2\tilde\beta_t(s)|=|p_2\beta_t(s)|$ for all
$s\in[0,\epsilon)$ and
$|\uparrow_{\gamma(t)}^{p_2}\uparrow_{\gamma(t)}^{\beta_t(s)}|=\frac\pi2$,
we have that (by the first variation formula)
$$\lim_{s\to0}\left|\uparrow_{\tilde\gamma(t)}^{\tilde p_2}
\uparrow_{\tilde\gamma(t)}^{\tilde\beta_t(s)}\right| =\frac{\pi}{2}\
(\text{i.e. } |\tilde\gamma'(t)\tilde X_{2m-1}(t)|=\frac\pi2)$$ (due
to Remark \ref{2.8}, one can also get this by Gauss's Lemma). Next we
will show that
$$\lim\limits_{s\to0}\left|\tilde X_j(t)\uparrow_{\tilde\gamma(t)}^{\tilde\beta_t(s)}\right|
=\frac{\pi}{2}\ (\text{i.e. } |\tilde X_j(t)\tilde
X_{2m-1}(t)|=\frac\pi2) \text{ for any } 1\leq j\leq 2m-2.\eqno{(3.6)}$$ Let $[p_{2j}^tq_{2j}^t]\subset
B_{M_2}(p_2,\frac{r_0}2)$ with
$\gamma(t)\in[p_{2j}^tq_{2j}^t]^\circ$ be a geodesic such that
$X_j(t)=\uparrow_{\gamma(t)}^{p_{2j}^t}$. From the choice of
$X_j(t)$, we know that $f_{[p_1p_2]}([p_{2j}^tq_{2j}^t])=[\tilde
p_{2j}^t\tilde q_{2j}^t]\subset\Bbb S^{2m-1}_{[p_1p_2],[p_2p_2']}$
with $|p_{2j}^tq_{2j}^t|=|\tilde p_{2j}^t\tilde q_{2j}^t|$, and that
$\tilde X_j(t)=\uparrow_{\tilde\gamma(t)}^{\tilde p_{2j}^t}$. Note
that $|\tilde p_{2j}^t\tilde \gamma(t)|=|p_{2j}^t\gamma(t)|, |\tilde
p_{2j}^t\tilde\beta_t(s)|\geq|p_{2j}^t\beta_t(s)|$ (by (2.8))  and
$|\uparrow_{\gamma(t)}^{p_{2j}^t}\uparrow_{\gamma(t)}^{\beta_t(s)}|=\frac\pi2$.
Then by the first variation formula we have that
$\lim\limits_{s\to0}|\uparrow_{\tilde\gamma(t)}^{\tilde p_{2j}^t}
\uparrow_{\tilde\gamma(t)}^{\tilde\beta_t(s)}|\geq\frac\pi2$, and
similarly $\lim\limits_{s\to0}|\uparrow_{\tilde\gamma(t)}^{\tilde
q_{2j}^t}\uparrow_{\tilde\gamma(t)}^{\tilde\beta_t(s)}|\geq\frac\pi2.$
Hence, (3.6) follows because $|\uparrow_{\tilde\gamma(t)}^{\tilde
p_{2j}^t}\uparrow_{\tilde\gamma(t)}^{\tilde q_{2j}^t}|=\pi$. On the
other hand, note that $|\tilde
\beta_t(s)\tilde\gamma(t)|\geq|\beta_t(s)\gamma(t)|$ (by (2.8)), and so
$|\tilde X_{2m-1}(t)|\geq|X_{2m-1}(t)|=1$.
(Moreover, it is easy to see that
$$\lim_{t\to0}|\tilde X_{2m-1}(t)|=1. \eqno{(3.7)}$$
In fact, if $m=1$ and if $(\rho,\theta)$ is the polar coordinates of
$M_2$ at $p_2$ in which $\gamma(t)$ has the coordinates $(t, 0)$ and
the metric $g_{M_2}=d\rho^2+G(\rho,\theta)d\theta^2$, then we have
that $\frac{|\tilde X_{2m-1}(t)|} {|X_{2m-1}(t)|}=\frac{\sin
t}{\sqrt{G(t,0)}}.$ So far, the proof of (and comments on) Claim 2
is finished.)

\vskip1mm

Now, we consider $f_{[p_1\gamma(t)]}$ and $\Bbb
S_{[p_1\gamma(t)]}^{2m}$. Let $\bar\beta_t(s)$ denote
$f_{[p_1\gamma(t)]}(\beta_t(s))$, which is a minimal geodesic in $\Bbb
S_{[p_1\gamma(t)]}^{2m}$ by (2.7). We can also define the corresponding
$\text{d}f_{[p_1\gamma(t)]}$, which is an isometrical embedding
(similar to Lemma \ref{2.6}). Hence, $\bar\beta_t'(0)$ is
perpendicular to $\tilde X_1(t),\cdots,\tilde X_{2m-2}(t)$ and
$\tilde\gamma'(t)$. On the other hand, note that $\tilde
X_1(t),\cdots,\tilde X_{2m-2}(t)$, $\frac{\tilde
X_{2m-1}(t)}{|\tilde X_{2m-1}(t)|}$ and $\tilde\gamma'(t)$ are
parallel and orthogonal along $\tilde\gamma(t)$ (on $\Bbb
S_{[p_1p_2]}^{2m}\subset\Bbb S_{p_1}^{2m+1}$). Then we can define an
orientable angle function $\theta(t)\in (-\pi,\pi]$ between
$\bar\beta_t'(0)$ and $\tilde X_{2m-1}(t)$. Note that
$$\theta(t)\neq 0, \pi \text{ for } t>0\eqno{(3.8)}$$
(otherwise it has to hold that $\Bbb S^{2m}_{[p_1\gamma(t)]}=\Bbb
S^{2m}_{[p_1p_2]}$, which contradicts (3.5)).

\noindent{\bf Claim 3}: {\it We have that
$$\cos\theta(t)=\frac{|X_{2m-1}(t)|}{|\tilde X_{2m-1}(t)|}
=\frac{1}{|\tilde X_{2m-1}(t)|}.\eqno{(3.9)}$$ As a corollary,
$\theta(t)$ is a smooth function (which implies that
$0<\theta(t)<\frac\pi2$ (or $-\frac\pi2<\theta(t)<0)$ and $|\tilde
X_{2m-1}(t)|>1$ for $t>0$ (see $(3.8)$).} By Lemma \ref{3.11}, we
first note that
$$\lim_{s\to0}\left|\uparrow_{\bar\beta_t(s)}^{\tilde\beta_t(s)}
\uparrow_{\bar\beta_t(s)}^{\tilde\gamma(t)}\right|=\frac\pi2.$$ Then
it is easy to see that
$$\cos\theta(t)=\lim_{s\to0}\frac{|\bar\beta_t(s)\tilde\gamma(t)|}
{|\tilde\beta_t(s)\tilde\gamma(t)|}=\frac{|X_{2m-1}(t)|}{|\tilde
X_{2m-1}(t)|}\ \text{ (i.e. (3.9) holds)}$$ (note that
$f_{[p_1\gamma(t)]}:\beta_t(s)|_{s\in[0,\epsilon)}\to\bar\beta_t(s)|_{s\in[0,\epsilon)}$
is an isometry). Due to (3.9), in order to prove that $\theta(t)$ is
a smooth function we need only to show that it is a continuous one.
We first observe that $\theta(t)\to 0$ as $t\to 0$ because
$\lim\limits_{t\to0}|\tilde X_{2m-1}(t)|=1$ (see (3.7)). Moreover,
note that $\theta(t+\Delta
t)=\theta(t)\pm\theta_{\tilde\gamma(t)}(\Delta t)$, where
$\theta_{\tilde\gamma(t)}(\Delta t)$ denotes the angle between
$\bar\beta_{t+\Delta t}'(0)$ and the vector (at
$\tilde\gamma(t+\Delta t)$) that is parallel to $\bar\beta_{t}'(0)$
along $\tilde\gamma$. Similar to $\theta(t)\to 0$ as
$t\to 0$, we have that $\theta_{\tilde\gamma(t)}(\Delta t)\to 0$
as $\Delta t\to 0$. That is, $\theta(t)$ is continuous with respect
to $t$. (Now, Claim 3 is verified.)

\vskip1mm

Based on Claim 3, we have the following important observation.

\noindent{\bf Claim 4}: {\it For any $q\in M_2$, $f_{p_1}(q)$ is a
closed $1$-dimensional smooth submanifold in $\Bbb S_{p_1}^{2m+1}$,
i.e. $f_{p_1}(q)$ consists of finite smooth circles which do not
intersect each other (and each of which does not intersect itself).}
(Of course, when the whole proof of the lemma has been finished, we
will know that $f_{p_1}(q)$ is just one or two great circles in
$\Bbb S_{p_1}^{2m+1}$.) According to Proposition \ref{3.1}, there
are $[qr], [rr'], [qr']\subset M_2$ with $|rq|,|rr'|<\frac{r_0}2$
such that the triangle formed by them does not bound a convex
spherical surface. Without loss of generality, we can assume that
$[rr']$ is just the $[p_2p_2']$. Then, due to (3.2), we have that
$f_{[p_1p_2]}(q)\not\in\Bbb S_{[p_1p_2],[p_2p_2']}^{2m-1}$.
According to (3.5), we have that $f_{[p_1\gamma(t)]}(q)\not\in\Bbb
S_{[p_1p_2],[p_2p_2']}^{2m-1}$ either. Note that $\gamma(t)$ belongs
to $B_{M_2}(q,r_0)$ for all $t$. It then follows from Claim 3 that
$$\alpha_q(t)|_{t\in[0,|p_2p_2'|]}\triangleq f_{[p_1\gamma(t)]}(q)|_{t\in[0,|p_2p_2'|]}
\text{ is a smooth curve.}\eqno{(3.10)}$$ Note that
$\alpha_q(t_1)\neq\alpha_q(t_2)$ for all $t_1\neq t_2$ by (3.5)
(note that $f_{[p_1\gamma(t)]}(q)\not\in\Bbb
S_{[p_1p_2],[p_2p_2']}^{2m-1}$),
so there is an interval $[a,b]\subseteq[0,|p_2p_2'|]$ such that the tangent vector
$$\alpha_q'(t)\neq0 \text{ for all } t\in[a,b].$$
By Lemma \ref{3.11}, $\alpha_q'(t)$ is perpendicular to $\Bbb
S^{2m}_{[p_1q]_t}$ for all $t\in[a,b]$, where $[p_1q]_t$ is the
minimal geodesic between $p_1$ and $q$ whose unit tangent vector at $p_1$ is
$\alpha_q(t)$. Moreover, note that
$$B_{\Bbb S^{2m}_{[p_1q]_t}}(\alpha_q(t),r_0)\cap
B_{\Bbb S^{2m}_{[p_1q]_{t'}}}(\alpha_q(t'),r_0)=\emptyset
\text{ for all } t\neq t'\in[a,b]$$ (this is due to that
$\alpha_q(t)\neq\alpha_q(t')$ and that $f_{[p_1\bar
p_2]}|_{B_{M_2}(\bar p_2,r_0)}$ is injective for any $[p_1\bar p_2]$
with $\bar p_2\in M_2$ (see Remark \ref{2.8})). Then
it is easy to see that the $r_0$-tubler neighborhood $U$ of $\alpha_q(t)|_{t\in[a,b]}$
satisfies $$U=\bigcup_{t\in [a,b]}B_{\Bbb S^{2m}_{[p_1q]_t}}(\alpha_q(t),r_0)
\text{ and } U\cap f_{p_1}(q)=\alpha_q(t)|_{t\in[a,b]}.\eqno{(3.11)}$$ On
the other hand, for any $\tilde q'\in f_{p_1}(q)$, there is
a minimal geodesic $[p_1p_2]'$ between $p_1$ and $p_2$
such that $f_{[p_1p_2]'}(q)=\tilde q'$. We also consider
$\Bbb S_{[p_1p_2]'}^{2m}$, the minimal geodesic $[p_1\gamma(t)]'$
between $p_1$ and $\gamma(t)$ whose unit tangent vector at $p_1$ is
$f_{[p_1p_2]'}(\gamma(t))$, and the
curve $\bar\alpha_q(t)|_{t\in[0,|p_2p_2'|]}\triangleq
f_{[p_1\gamma(t)]'}(q)|_{t\in[0,|p_2p_2'|]}$ (with
$\bar\alpha(0)=\tilde q'$). Note that
$$\bar\alpha_q(t)|_{t\in[0,|p_2p_2'|]}\text{ is identical to
$\alpha_q(t)|_{t\in[0,|p_2p_2'|]}$ up to an isometry of $\Bbb
S_{p_1}^{2m+1}$}.\eqno{(3.12)}$$ This together with (3.11) implies
that $f_{p_1}(q)$ is a closed 1-dimensional submanifold in $\Bbb
S_{p_1}^{2m+1}$. (The proof of Claim 4 is done.)

\vskip1mm

Based on Claim 4, we can draw an almost immediate conclusion.

\noindent{\bf Claim 5}: {\it $f_{p_1}:M_2\to \Bbb S^{2m+1}_{p_1}$ is
surjective.} Let $\cal S$ be a component of $f_{p_1}(q)$ for the $q$
in the proof of Claim 4,
which is a smooth circle in $\Bbb S^{2m+1}_{p_1}$. For any $\tilde
z\in \Bbb S^{2m+1}_{p_1}$, there is a $\tilde q\in\cal S$ such that
$|\tilde z\tilde q|=\min\{|\tilde z\tilde q'||\tilde q'\in\cal S\}$.
Note that $[\tilde q\tilde z]$ is perpendicular to $\cal S$ at
$\tilde q$. On the other hand, there is a $[p_1q]$ such
that $\tilde q=\uparrow_{p_1}^q$, and $\cal S$ is perpendicular to
$\Bbb S^{2m}_{[p_1q]}$ at $\tilde q$ (by Lemma \ref{3.11}). It then
follows that $\tilde z$ belongs to $\Bbb S^{2m}_{[p_1q]}\subset f_{p_1}(M_2)$
(see (2.10.1)).

\vskip1mm

We still consider the $\cal S$ in the proof of Claim 5. Let
$s\in[0,\ell]$ be the arc-length parameter of $\cal S$, where $\ell$ is the
perimeter of $\cal S$. It follows from the proof Claim 5 that
$$\Bbb S_{p_1}^{2m+1}=\bigcup_{s\in[0,\ell]}\Bbb S^{2m}_{[p_1q]_s},$$
where $[p_1q]_s$ is the minimal geodesic between $p_1$ and $q$ whose unit
tangent vector at $p_1$ is ${\cal S}(s)$. Note that there is a natural
isometry  (ref. Remark 2.8)
$$h_{s,s'}\triangleq \exp_{{\cal S}(s')}\circ\text{ d}f_{[p_1q]_{s'}}
\circ\text{d}f_{[p_1q]_s}^{-1}\circ\exp_{{\cal
S}(s)}^{-1}:\Bbb S^{2m}_{[p_1q]_s}\to\Bbb
S^{2m}_{[p_1q]_{s'}}.\eqno{(3.13)}$$ By Lemma \ref{2.4} (and (2.5)), we
observe that $h_{s,s'}\to h_{s,s_0}$ if $s'\to s_0$, and thus
$h_{0,s}(\tilde x)|_{s\in[0,\ell]}$ is continuous with respect to
$s$ for any $\tilde x\in \Bbb S^{2m}_{[p_1q]_0}$. Observe that
$$h_{0,s}(\tilde x)\neq h_{0,s'}(\tilde x) \text{ for any }0\leq s\neq
s'<\ell.\eqno{(3.14)}$$ (Otherwise,
$h_{0,s}([{\cal S}(0)\tilde x])$ and $h_{0,s'}([{\cal
S}(0)\tilde x])$ are two minimal geodesics starting from $h_{0,s}(\tilde x)$
in $\Bbb S^{2m}_{[p_1x]_s}$, where $x=f_{p_1}^{-1}(\tilde x)$
and $[p_1x]_s$ is the minimal geodesic between $p_1$ and $x$
whose unit tangent vector at $p_1$ is $h_{0,s}(\tilde x)$. However, note that
$f_{p_1}^{-1}(h_{0,s}([{\cal S}(0)\tilde
x]))=f_{p_1}^{-1}(h_{0,s'}([{\cal S}(0)\tilde
x]))=f_{p_1}^{-1}([{\cal S}(0)\tilde x])$. We know that this is impossible by
Remark \ref{2.8} when we consider $f_{[p_1x]_s}$.) It follows that
$$h_{0,s}(\tilde x)|_{s\in[0,\ell]}\text{ is a component of } f_{p_1}(x).\eqno{(3.15)}$$

On the other hand, note that there is a neighborhood $V$ of $q$ in
$B_{M_2}(p_2,\frac{r_0}{2})$ such that $f_{[p_1p_2]}(V)\cap\Bbb
S_{[p_1p_2],[p_2p_2']}^{2m-1}=\emptyset$ (because
$f_{[p_1p_2]}(q)\not\in\Bbb S_{[p_1p_2],[p_2p_2']}^{2m-1}$).
Similar to $\alpha_{q}(t)$ and $\bar\alpha_{q}(t)$, we can define
$\alpha_{v}(t)|_{t\in[0,|p_2p_2'|]}\triangleq
f_{[p_1\gamma(t)]}(v)|_{t\in[0,|p_2p_2'|]}$ and
$\bar\alpha_{v}(t)|_{t\in[0,|p_2p_2'|]}\triangleq
f_{[p_1\gamma(t)]'}(v)|_{t\in[0,|p_2p_2'|]}$ for any $v\in V$. In
fact, for any $v_1, v_2\in V$, we have a strengthened version of
(3.12) (note that the isometry of $\Bbb S_{p_1}^{2m+1}$ mentioned in
(3.12) restricted to $\Bbb S_{[p_1p_2]}^{2m}$ is actually an
isometry from $\Bbb S_{[p_1p_2]}^{2m}$ to $\Bbb S_{[p_1p_2]'}^{2m}$
like (3.13)):
$$|\alpha_{v_1}(t_1)\alpha_{v_2}(t_2)|=|\bar\alpha_{v_1}(t_1)\bar\alpha_{v_2}(t_2)|
\text{ for all } t_1,t_2\in[0,|p_2p_2'|].\eqno{(3.16)}$$ (This is a
very important observation to the whole proof.)

Based on (3.16), we can conclude that: if $0<\Delta s<\delta$ for a small
$\delta$, then the map
$$h_{0,s}(\tilde x)|_{s\in[0,\Delta s]}\to h_{0,s}(\tilde x)|_{s\in[\Delta s,2\Delta s]}
\text{ defined by } h_{0,s}(\tilde x)\mapsto h_{s,s+\Delta
s}(h_{0,s}(\tilde x)) \text{ is an isometry}\eqno{(3.17)}$$ for all
$\tilde x\in \Bbb S^{2m}_{[p_1q]_0}$ (note that $s$ is the arc-length parameter of $\cal S$). Let ${\cal S}^*$ denote the circle
$h_{0,s}(\tilde x)|_{s\in[0,\ell]}$ (see (3.15) and Claim 4), and
let $s^*$ be its arc-length parameter. Note that we can assume that ${\cal
S}^*(0)=h_{0,0}(\tilde x)=\tilde x$, and ${\cal S}^*(\Delta
s^*)=h_{0,\Delta s}(\tilde x)$ for some $\Delta s^*>0$. It then
follows from (3.17) and (3.14) that $$h_{0,s}(\tilde
x)|_{s\in[\Delta s,2\Delta s]} \text{ is the arc } {\cal
S}^*(s^*)|_{s^*\in[\Delta s^*,2\Delta s^*]} \text{ with }
h_{0,2\Delta s}(\tilde x)={\cal S}^*(2\Delta s^*)$$ (NOT the arc
${\cal S}^*(s^*)|_{s^*\in[0,\Delta s^*]}$). Similarly, we have that
$h_{0,k\Delta s}(\tilde x)={\cal S}^*(k\Delta s^*)$ for any $k\in
\Bbb N^+$. Due to the arbitrariness of $\Delta s$, this implies that
$$\frac{\Delta s^*}{\Delta
s}=\frac{\ell^*}{\ell} \text{ and } h_{0,s}(\tilde x)={\cal
S}^*(s^*) \text{ with } \frac{s^*}{s}=\frac{\ell^*}{\ell} \text{ for
any } s\in[0,\ell],\eqno{(3.18)}$$ where $\ell^*$ is the perimeter
of ${\cal S}^*$. Consequently, if $h_{0,s_0}(\tilde x)={\cal
S}^*(s^*_0)$, then $$h_{s_0,s_0+s}({\cal S}^*(s^*_0))={\cal
S}^*(s^*_0+\frac{\ell^*}{\ell}s).\eqno{(3.19)}$$ An important fact
is that if ${\cal S}^*(s^*_0)\in \Bbb S^{2m}_{[p_1q]_0}$ for some
$0<s^*_0<\ell^*$, we also have that
$$h_{0,s}({\cal S}^*(s^*_0))={\cal S}^*(s^*_0+\frac{\ell^*}{\ell}s).\eqno{(3.20)}$$
If this is not true, then, by (3.18), it has to hold that
$h_{0,s}({\cal S}^*(s^*_0))={\cal S}^*(s^*_0-\frac{\ell^*}{\ell}s)$,
which implies that $h_{0,\frac{s_0}{2}}(\tilde
x)=h_{0,\frac{s_0}{2}}({\cal S}^*(s^*_0))={\cal
S}^*(\frac{s^*_0}{2})$. Since $h_{0,\frac{s_0}{2}}$ is an isometry,
we have that $\tilde x={\cal S}^*(s^*_0)=h_{0,s_0}(\tilde x)$, which
contradicts (3.14). On
the other hand, by Claims 5 and 4 we know that $$\Bbb
S_{p_1}^{2m+1}=\bigcup_{p_2\in M_2}f_{p_1}(p_2),$$ where each
$f_{p_1}(p_2)$ consists of finite components similar to ${\cal S}^*$.
Moreover, by the proof of Claim
5, there exist a $\tilde z_0\in \Bbb S^{2m}_{[p_1q]_0}$ and
$s_{\tilde z_0}$ such that $h_{0,s_{\tilde z_0}}(\tilde z_0)=\tilde
z$ for any $\tilde z\in\Bbb S_{p_1}^{2m+1}$.
Therefore, based on $\cal S$ and due to (3.19-20), we can define an $S^1$-action on $\Bbb
S^{2m+1}_{p_1}$
$$h:{\cal S}\times \Bbb S^{2m+1}_{p_1}\rightarrow \Bbb S^{2m+1}_{p_1}
\text{ defined by } h(s,\tilde z)=h_{0,s_{\tilde z_0}+s}(\tilde
z_0),\eqno{(3.21)}$$ where $\tilde z_0\in \Bbb S^{2m}_{[p_1q]_0}$ and
$h_{0,s_{\tilde z_0}}(\tilde z_0)=\tilde z$.

\noindent{\bf Claim 6}: {\it Through $h$, $\cal S$ acts on $\Bbb
S^{2m+1}_{p_1}$ freely and isometrically.} By (3.14), $\cal S$ acts
on $\Bbb S^{2m+1}_{p_1}$ freely. It then suffices to show that each
$h(s,\cdot)$ is a local isometry (note that $h(s,\cdot)$ is a 1-1
map). For any $\tilde z\in\Bbb S_{p_1}^{2m+1}$ and $z\triangleq
f_{p_1}^{-1}(\tilde z)$, we let ${\cal S}_{\tilde z}$ be the
component of $f_{p_1}(z)$ containing $\tilde z$. We first assume
that $\tilde z$ is sufficiently close to some ${\cal S}(s_0)$, and
will prove that
$$|h(s,{\cal S}(s_0))h(s,\tilde z)|=|{\cal S}(s_0)\tilde z|.\eqno{(3.22)}$$
Note that $h(s,{\cal S}(s_0))={\cal S}(s_0+s)$ and
$h(s,\tilde z)=h_{0,s_{\tilde z_0}+s}(\tilde z_0)$, and that
$\tilde z$ and $h(s,\tilde z)$ are the unique points in
${\cal S}_{\tilde z}$ such that $$|{\cal S}(s_{\tilde
z_0})\tilde z|=|{\cal S}(s_{\tilde
z_0}+s)h(s,\tilde z)|=|{\cal S}(0)\tilde
z_0|\ (\leq|{\cal S}(s_0)\tilde z|)$$ (see Remark \ref{2.8}).
Of course, this implies that $|{\cal S}(s_{\tilde
z_0}){\cal S}(s_0)|$ is sufficiently small, so is $|s_{\tilde z_0}-s_0|$.
It then is not hard to see that (3.22) follows from (3.16).
Now we let $\tilde z$ be an arbitrary point in $\Bbb S_{p_1}^{2m+1}$,
and we need only to show that
$$|h(s,\tilde z')h(s,\tilde z)|=|\tilde z'\tilde z|\eqno{(3.23)}$$
for any $\tilde z'$ in $\Bbb S_{p_1}^{2m+1}$ sufficiently close to
$\tilde z$. Let $s^*$ (resp. $s^{**}$) be the arc-length parameter of
${\cal S}_z$ (resp. ${\cal S}_{z'}$),  which is increasing with respect to $s$,
such that $h(s,\tilde z)={\cal S}_{\tilde z}(s^*)$
with ${\cal S}_{\tilde z}(0)=\tilde z$
(resp. $h(s,\tilde z')={\cal S}_{\tilde z'}(s^{**})$ with ${\cal S}_{\tilde z'}(0)=\tilde z'$).
If we replace $\cal S$ with ${\cal S}_{\tilde z}$, we can similarly define
$\bar h(s^*,\cdot)$ such that $\bar h(s^*,\tilde z)={\cal S}_{\tilde z}(s^*)$.
Similarly, we have that  $\bar h(s^*,\tilde z')={\cal S}_{\tilde z'}({s^{**}}')$
with ${\cal S}_{\tilde z'}({s^{**}}')|_{{s^{**}}'=0}=\tilde z'$,
where ${s^{**}}'$ is another arc-length parameter of ${\cal S}_{\tilde z'}$;
and by (3.22) we have that
$$|\bar h(s^*,\tilde z)\bar h(s^*,\tilde z')=|\tilde z\tilde z'|.\eqno{(3.24)}$$
Note that ${s^{**}}'=s^{**}$ or ${s^{**}}'=-s^{**}$, and it will
suffice to show that the latter case does not occur. In fact, if
${s^{**}}'=-s^{**}$, then (3.24) implies that $|h(s,\tilde z_1)h(s,
\tilde z')|$ will change as $s$ changes, where $\tilde z_1$ is a
point in ${\cal S}_{\tilde z}$ such that  $\tilde z_1$ and $\tilde
z'$ lie in some $\Bbb S^{2m}_{[p_1q]_{s_1}}$. This is impossible
because $|h(s,\tilde z_1)h(s, \tilde z')|=|h_{s_1,s_1+s}(\tilde
z_1)h_{s_1,s_1+s}(\tilde z')|$ and $h_{s_1,s_1+s}$ is an isometry.
Note that the proof of Claim 6 is completed now.

\vskip2mm

By Claim 6, we know that $\Bbb S^{2m+1}_{p_1}/{\cal S}=\Bbb{CP}^m.$
And it is easy to see that
$$\Bbb S^{2m+1}_{p_1}/{\cal S}(=\Bbb{CP}^m)\to M_2
\text{ defined by } \{f_{p_1}(p_2)\}\mapsto p_2 \eqno{(3.25)}$$ is a
locally isometrical map (see (2.2) and Lemma 2.3), i.e. it is a Riemannian covering map.
Therefore, according to Synge's theorem ([CE]), $M_2$
is isometric to ${\Bbb C\Bbb P}^{m}$ or ${\Bbb C\Bbb P}^{m}/\Bbb
Z_2$ with the canonical metric, and that $M_2\stackrel{\rm iso}{\cong}\Bbb
{CP}^{m}/\Bbb Z_2$ occurs only when $m$ is odd (see Lemma A.1
in Appendix). (The long proof of Lemma \ref{3.10} is completed now.)

\hfill$\Box$

\vskip2mm

Based on Lemma \ref{3.10} and its proof, we give the following
two important facts. In the following, we assume that $n_i=2m_i$ for
$i=1$ and 2 (see Lemma \ref{3.8}).

\begin{lemma}\label{3.12}
For any $p\in M$, there exists a $[p_1p_2]$ with $p_i\in M_i$ such
that $p\in[p_1p_2]$.
\end{lemma}

\noindent{\bf Proof.} We select $p_1\in M_1$ such that
$|pp_1|=\min\{|pp_1'||p_1'\in M_1\}$. Note that for any $[p_1p]$, we
have that $\uparrow_{p_1}^p\in(\Sigma_{p_1}M_1)^{=\frac{\pi}{2}}$
($=\Bbb S_{p_1}^{2m_2+1}$). By Claim 5 in the proof of Lemma
\ref{3.10}, there is a $[p_1p_2]$ with $p_2\in M_2$ such that
$$\uparrow_{p_1}^p=\uparrow_{p_1}^{p_2}.$$
{\bf Claim}: {\it $|pp_1|\leq\frac\pi2$, and thus $p\in [p_1p_2]$.}
If $|pp_1|>\frac\pi2$, then $p_2\in [p_1p]$, and so $|pp_2|\leq\frac\pi2$ (because
$|p_1p|\leq \pi$). Note that
$|pp_2|=\min\{|pp_2'||p_2'\in M_2\}$ because $|p_1p_2'|=\frac\pi2$
for all $p_2'\in M_2$. Similarly, if
$$f_{p_2}(M_1)=(\Sigma_{p_2}M_2)^{=\frac{\pi}{2}},\eqno{(3.26)}$$
then we can find a $[p_2\bar p_1]$ with $\bar p_1\in M_1$ such that $p\in [p_2\bar
p_1]$, which contradicts the choice of $p_1$. Hence, we need only to prove
(3.26). If $n_1=2m_1\geq2$, then (3.26) automatically holds (similar to
Claim 5 in the proof of Lemma \ref{3.10}). If $m_1=0$ (i.e.
$M_1=\{p_1\}$), there is a natural map
$$\tau: f_{p_1}(p_2)\to (\Sigma_{p_2}M_2)^{=\frac{\pi}{2}}
\text{ defined by }
\uparrow_{p_1}^{p_2}\mapsto\uparrow_{p_2}^{p_1},$$ where
$\uparrow_{p_1}^{p_2}$ and $\uparrow_{p_2}^{p_1}$ are the directions
of any given $[p_1p_2]$. Obviously, $\tau$ is injective and
continuous. On the other hand, note that $(\Sigma_{p_2}M_2)^{=\frac{\pi}{2}}$
is a circle (in this case $M_2$ is of codimension 2), and each
component of $f_{p_1}(p_2)$ is a circle (see Claim 4 in the
proof of Lemma \ref{3.10}). Hence, $\tau$ restricted
to each component of $f_{p_1}(p_2)$ is surjective, and thus (3.26) follows
(and $f_{p_1}(p_2)$ contains only one component).
\hfill$\Box$

\begin{lemma}\label{3.13}
If $M_2$ is isometric to ${\Bbb C\Bbb P}^{m_2}$ (here $m_2$ may be
$0$), then

\noindent{\rm (3.13.1)} $M_1$ is isometric to ${\Bbb C\Bbb
P}^{m_1}$.

\noindent{\rm (3.13.2)} For any $p\in M$, $\{p\}^{=\frac\pi2}$ is
totally geodesic and of codimension $2$ in $M$, and so it is
isometric to $\Bbb{CP}^{\frac n2-1}$.
\end{lemma}

\noindent{\bf Proof.}  (3.13.1) We need only to consider the case
$m_1>0$. We consider the following natural map
$$\bar\tau: f_{p_2}(p_1)\to f_{p_1}(p_2)
\text{ defined by }
\uparrow_{p_2}^{p_1}\mapsto\uparrow_{p_1}^{p_2}$$ (similar to the
$\tau$ in the proof of Lemma \ref{3.12}). Since $M_2$ is isometric
to ${\Bbb C\Bbb P}^{m_2}$, from the end of the proof of Lemma
\ref{3.10} (resp. Lemma \ref{3.12}), we can conclude that
$f_{p_1}(p_2)$ is just a circle when $m_2>0$ (resp. $m_2=0$). Then
by the arguments in the end of the proof of Lemma \ref{3.12},
$f_{p_2}(p_1)$ contains only one component (so it is a circle
too). This together with Lemma \ref{3.10} (and the end of its proof)
implies that $M_1\stackrel{\text{iso}}{\cong}{\Bbb C\Bbb P}^{m_1}$.

\vskip1mm

(3.13.2) We first supply the proof for any given point $p\in M_i$,
say $M_2$. Since $M_2\stackrel{\text{iso}}{\cong}\Bbb{CP}^{m_2}$,
$N_2\triangleq\{p\}^{=\frac\pi2}\cap M_2$ is isometric to
$\Bbb{CP}^{m_2-1}$ and is totally geodesic in $M_2$ (note that we
need only to consider the case ``$m_2>0$''). Then, by Lemma
\ref{3.12} and \ref{2.3}, it is not hard to see that
$$\{p\}^{=\frac\pi2}=\{q\in M|q \text{ belongs to some $[p_1p_2]$ with
$p_1\in M_1$ and $p_2\in N_2$}\}$$ and that
$$\{p\}^{=\frac\pi2}=\{p\}^{\geq\frac\pi2}.\eqno{(3.27)}$$
Moreover, we have that $\lambda_{pq}=+\infty$ for any $q\in
\{p\}^{=\frac\pi2}$. Hence, by (3.13.1) we can get that
$\{p\}^{=\frac\pi2}$ is isometric to $\Bbb{CP}^{\frac n2-1}$ once we
have proved that it is totally geodesic (note that it is closed in
$M$) and of codimension $2$ in $M$.

Note that $\{p\}^{\geq\frac\pi2}$ is convex in $M$ by (i) of Theorem
1.1. Then (3.27) and ``$\lambda_{pq}=+\infty\ \forall\ q\in
\{p\}^{=\frac\pi2}$'' imply that $\dim(\{p\}^{=\frac\pi2})\leq n-2$
(ref. [RW]\footnote{In [RW], it has been proved that: {\it Let $A_1$
and $A_2$ be two convex subsets in an $n$-dimensional Alexandrov
space with curvature $\geq 1$. If $|a_1a_2|=\frac\pi2$ for any
$a_i\in A_i$, then $\dim(A_1)+\dim(A_2)\leq n-1$; and if equality
holds, then $\lambda_{a_1a_2}<+\infty$ for all $a_i\in A_i^\circ$
(where $X^\circ$ denotes the interior part of $X$).} }). On the
other hand, note that both $M_1$ and $N_2$ are totally geodesic in
$\{p\}^{=\frac\pi2}$, and that $|p_1p_2|=\frac\pi2$ and
$\lambda_{p_1p_2}=+\infty$ for all $p_1\in M_1$ and $p_2\in N_2$.
Similarly, it is implied that $\dim(\{p\}^{=\frac\pi2})\geq n-2$. It
then follows that $\dim(\{p\}^{=\frac\pi2})=n-2.$

Next we will prove that $\{p\}^{=\frac\pi2}$ has empty boundary,
which implies that $\{p\}^{=\frac\pi2}$ is totally geodesic (because
it is convex) in $M$. Since any $q\in\{p\}^{=\frac\pi2}$ lies in a
$[p_1p_2]$ with $p_1\in M_1$ and $p_2\in N_2$, it suffices to
show that both $M_1$ and $N_2$ consist of interior  points\footnote{We
know that, in an Alexandrov space $A$ with curvature bounded below,
any minimal geodesic between two interior  points belongs to $A^\circ$
(ref. [BGP]).} of $\{p\}^{=\frac\pi2}$. Let $p_2$ be an arbitrary
point in $N_2$. By (3.26), $(\Sigma_{p_2}M_2)^{=\frac{\pi}{2}}$ belongs to
$\Sigma_{p_2}\{p\}^{=\frac\pi2}$. Furthermore, by Lemma \ref{2.3}
it is easy to see that
$$\Sigma_{p_2}\{p\}^{=\frac\pi2}=(\Sigma_{p_2}M_2)^{=\frac{\pi}{2}}*\Sigma_{p_2}N_2=\Bbb S^{n-3}$$
(note that $N_2$ is totally geodesic in $\{p\}^{=\frac\pi2}$).
It follows that $p_2$ $(\in N_2)$ is an interior
point of $\{p\}^{=\frac\pi2}$. Let $p_1$ be an arbitrary point in
$M_1$. From (3.25), it is easy to see that $f_{p_1}(N_2)$ is an
$\Bbb S^{2m_2-1}$ in $\Bbb S_{p_1}^{2m_2+1}$. Then similarly, we can
get that
$\Sigma_{p_1}\{p\}^{=\frac\pi2}=f_{p_1}(N_2)*\Sigma_{p_1}M_1=\Bbb
S^{n-3}$, i.e. $p_1$ is also an interior point of $\{p\}^{=\frac\pi2}$.

So far we have given the proof for any point $p\in M_2$. Now let $p$
be an arbitrary point in $M$. By Lemma \ref{3.12}, there is a
$[p_1p_2]$ with $p_i\in M_i$ such that $p\in [p_1p_2]$. Since
$M_i\stackrel{\text{\rm iso}}{\cong}{\Bbb{CP}^{m_i}}$ with $m_1$ or
$m_2>0$, say $m_2>0$, we can select $\bar p_2\in M_2$ such that
$|\bar p_2p_2|=\frac\pi2$. We have proved that $\{\bar
p_2\}^{=\frac\pi2}$ is isometric to $\Bbb{CP}^{\frac n2-1}$, and it
is easy to see that $p\in \{\bar p_2\}^{=\frac\pi2}$ by Lemma
\ref{2.3}. Hence, (3.13.2) follows if we replace $M_1$ and $M_2$
with $\bar M_1\triangleq\{\bar p_2\}$ and $\bar M_2\triangleq\{\bar
p_2\}^{=\frac\pi2}$ respectively. \hfill$\Box$

\vskip2mm

To the whole proof of the Main Theorem, the most difficult parts are
to prove that $M_i$ ($i=1,2$) are both isometric to ${\Bbb C\Bbb
P}^{m_i}$ or ${\Bbb C\Bbb P}^{m_i}/\Bbb Z_2$, and that any $p\in M$
lies in a $[p_1p_2]$ with $p_i\in M_i$. We would like to point out
that once these are established, we can find an argument in [GG1] to
prove that $M$ isometric to ${\Bbb C\Bbb P}^{\frac n2}$ or ${\Bbb
C\Bbb P}^{\frac n2}/\Bbb Z_2$, i.e. the following lemma holds. For
the convenience of readers and the completeness of the present
paper, we will supply a detailed proof for it.

\begin{lemma}\label{3.14}
If $\lambda_{p_1p_2}=+\infty$ for all $p_i\in M_i$, then
$M_i\stackrel{\text{\rm iso}}{\cong}{\Bbb C\Bbb P}^{m_i}$ and
$M\stackrel{\text{\rm iso}}{\cong}{\Bbb C\Bbb P}^{\frac n2}$, or
$M_i\stackrel{\text{\rm iso}}{\cong}{\Bbb C\Bbb P}^{m_i}/\Bbb Z_2$ and
$M\stackrel{\text{\rm iso}}{\cong}{\Bbb C\Bbb P}^{\frac n2}/\Bbb Z_2$
(only when $m_i$ and $\frac{n}2$ are odd) with
canonical metrics.
\end{lemma}

\noindent{\bf Proof.} By Lemma \ref{3.10} and \ref{3.13}, $M_i$ is
isometric to ${\Bbb C\Bbb P}^{m_i}$ or ${\Bbb C\Bbb P}^{m_i}/\Bbb
Z_2$ ($i=1, 2$) with the canonical metric, and that
$M_i\stackrel{\text{\rm iso}}{\cong}{\Bbb C\Bbb P}^{m_i}/\Bbb Z_2$
occurs only when $m_1$ and $m_2$ are odd. Then we can divide the
proof into the following two cases.

\vskip1mm

\noindent{\bf Case 1}. $M_i\stackrel{\rm iso}{\cong}{\Bbb C\Bbb
P}^{m_i}$. In this case, we will prove that $M$ is isometric to
${\Bbb C\Bbb P}^{\frac n2}$.

According to (3.13.2), we can assume that $n_1=0$ (i.e.
$M_1=\{p_1\}$) and $n_2=n-2$, and thus
$M_2\stackrel{\text{iso}}{\cong}\Bbb{CP}^{\frac n2-1}$. Let
$\nu:M_2\hookrightarrow\Bbb{CP}^{\frac n2}$ be an isometrical
embedding whose image is denoted by $\hat M_2$, and let $\hat p_1$
be the point in $\Bbb{CP}^{\frac n2}$ such that $d(\hat p_1,\hat
M_2)=\frac\pi2$ (note that $|\hat p_1\hat p_2|=\frac\pi2$ for any
$\hat p_2\in \hat M_2$). Due to Claim 5 and (3.25) in the proof of
Lemma \ref{3.10}, $\Sigma_{p_1}M$ (resp. $\Sigma_{\hat
p_1}\Bbb{CP}^{\frac n2}$) admits an isometrical and free
$S^1$-action such that each $S^1$-orbit is some
$\Uparrow_{p_1}^{p_2}$ (resp. $\Uparrow_{\hat p_1}^{\hat p_2}$) and
$(\Sigma_{p_1}M)/S^1=M_2$ (resp. $(\Sigma_{\hat p_1}\Bbb{CP}^{\frac
n2})/S^1=\hat M_2$). Hence, there is a natural isometrical map
$$\text{\bf i}_*:\Sigma_{p_1}M\to\Sigma_{\hat p_1}\Bbb{CP}^{\frac n2}\ (=\Bbb
S^{n-1})$$ such that $\text{\bf i}_*(\Uparrow_{p_1}^{p_2})=\Uparrow_{p_1}^{\hat
p_2}$ with $\hat p_2=\nu(p_2)$ for any $p_2\in M_2$. Furthermore,
due to Lemma \ref{3.12}, $\text{\bf i}_*$ induces a natural 1-1 map
$$\text{\bf i}:M\to \Bbb{CP}^{\frac n2}$$ with
$\text{\bf i}(p_1)=\hat p_1$, $\text{\bf i}|_{M_2}=\nu$ and
$\text{\bf i}([p_1p_2])=[\hat p_1\hat p_2]$ for any $[p_1p_2]$ with $p_2\in M_2$
such that $\uparrow_{\hat p_1}^{\hat p_2}=\text{\bf
i}_*(\uparrow_{p_1}^{p_2})$ and $\text{\bf i}|_{[p_1p_2]}$ is an isometry.
Let $x$ and $y$ be any two points in $M$. We need to show that
$$|\text{\bf i}(x)\text{\bf i}(y)|=|xy|.\eqno{(3.28)}$$
By Lemma \ref{3.12}, we can select $[p_1p_x]$ and $[p_1p_y]$ with
$p_x,p_y\in M_2$ such that $x\in [p_1p_x]$ and $y\in [p_1p_y]$. In
the following, $\hat p$ always denotes $\text{\bf i}(p)$ for any
$p\in M$. Note that $\hat x\in [\hat p_1\hat p_x]$ and $\hat y\in
[\hat p_1\hat p_y]$ with $|\hat x\hat p_1|=|xp_1|$ and $|\hat y\hat
p_1|=|yp_1|$. Since $M_2\stackrel{\text{iso}}{\cong}\Bbb{CP}^{\frac
n2-1}$, there is a $[p_xp_2]\subset M_2$ with $|p_xp_2|=\frac\pi2$
such that $p_y\in[p_xp_2]$. By (iii) of Theorem 1.1, there are
$[p_1p_2]$ and another minimal geodesic $[p_1p_x]'$ between $p_1$
and $p_x$ such that the triangle formed by $[p_1p_2]$, $[p_1p_x]'$
and $[p_xp_2]$ bounds a convex spherical surface which contains
$[p_1p_y]$ (ref. [GM]). In this surface, there is a $[p_2y']$ with $y'\in
[p_1p_x]'$ such that $y\in [p_2y']$. And an important point is that,
based on Lemma \ref{2.4}, it is not hard to see that $\hat y$
belongs to $[\hat p_2\hat{y'}]$ with $|\hat y\hat{y'}|=|yy'|$
(note that $\hat{y'}\in \text{\bf
i}([p_1p_x]')$).

Note that $x,y'\in \{p_2\}^{=\frac\pi2}$, so $[p_2y']$ is
perpendicular to $\{p_2\}^{=\frac\pi2}$ at $y'$ (note that
$\{p_2\}^{=\frac\pi2}$ is totally geodesic in $M$ by (3.13.2)). Then
by Lemma \ref{2.3}, it is easy to see that
$$\cos|xy|=\cos|yy'|\cos|xy'|.$$
On the other hand, similarly, it is not hard to
see that
$$\cos|\hat x\hat y|=\cos|\hat y\hat{y'}|\cos|\hat x\hat{y'}|$$
(with $|\hat y\hat{y'}|=|yy'|$). Hence, in order to see (3.28), it
suffices to show that $$|\hat x\hat{y'}|=|xy'|.\eqno{(3.29)}$$ By
the definition of $\text{\bf i}$, it is easy to see that $\text{\bf
i}({\{p_2\}^{=\frac\pi2}})=\{\hat p_2\}^{=\frac\pi2}$. Moreover, by
(3.13.2) we know that $\{p_2\}^{=\frac\pi2}$ is isometric to
$\Bbb{CP}^{\frac n2-1}$, which together with (3.25) implies that
$\text{\bf i}|_{{\{p_2\}^{=\frac\pi2}}}$ is an isometry. Hence,
(3.29) follows (because $x,y'\in \{p_2\}^{=\frac\pi2}$).

\vskip1mm

\noindent {\bf Case 2}. $M_i\stackrel{\rm iso}{\cong}{\Bbb C\Bbb
P}^{m_i}/\Bbb Z_2$. In this case, $\frac n2$ is odd because both
$m_1$ and $m_2$ are odd, and we will prove that $M$ is isometric to
${\Bbb C\Bbb P}^{\frac n2}/\Bbb Z_2$.

Note that $M_i$ $(i=1,2)$ can be embedded isometrically into
$\Bbb{CP}^{\frac n2}/\Bbb Z_2$ with $|\hat p_1\hat
p_2|=\frac\pi2$ for any $\hat p_i\in \hat M_i$
(see (A.1) in Appendix), where $\hat M_i$ (resp. $\hat p_i$) denotes the embedding
image of $M_i$ (resp. any given point $p_i\in M_i$).
Similar to the $\text{\bf i}_*$ in Case 1, for a
fixed point $p_{2,0}\in M_2$, there is an isometrical map
$$j_*:(\Sigma_{p_{2,0}}M_2)^{=\frac{\pi}{2}}\ (=\Bbb
S_{p_{2,0}}^{2m_1+1})\to(\Sigma_{\hat p_{2,0}}\hat
M_2)^{=\frac{\pi}{2}}\ (=\Bbb S_{\hat p_{2,0}}^{2m_1+1})$$ with
$j_*(\Uparrow_{p_{2,0}}^{p_1})=\Uparrow_{\hat p_{2,0}}^{\hat p_1}$
for any $p_1\in M_1$. Note that $j_*|_{\Uparrow_{p_{2,0}}^{p_1}}$
induces a natural homeomorphism $\bar
j_{*,p_1}:\Uparrow_{p_1}^{p_{2,0}}\to \Uparrow_{\hat p_1}^{\hat
p_{2,0}}$ which maps the unit tangent vector at $p_1$ of a $[p_1p_{2,0}]$ to
that at $\hat p_1$ of the $[\hat p_1\hat p_{2,0}]$ with $\uparrow_{\hat
p_{2,0}}^{\hat p_1}=j_*(\uparrow_{p_{2,0}}^{p_1})$. It
is not hard to see that, for any $p_1\in M_1$, there is a unique homeomorphism (ref. the
$\text{\bf i}_*$ in Case 1)
$$i_{p_1*}:(\Sigma_{p_1}M_1)^{=\frac{\pi}{2}}\ (=\Bbb
S_{p_1}^{2m_2+1})\to(\Sigma_{\hat p_1}\hat M_1)^{=\frac{\pi}{2}}\ (=\Bbb
S_{\hat p_1}^{2m_2+1})$$ such that
$i_{p_1*}|_{\Uparrow_{p_1}^{p_{2,0}}}=\bar j_{*,p_1}$,
$i_{p_1*}(\Uparrow_{p_1}^{p_{2}})=\Uparrow_{\hat p_1}^{\hat p_{2}}$ for any
$p_2\in M_2$, and $|i_{p_1*}(\uparrow_{p_1}^{p_{2}})i_{p_1*}(\uparrow_{p_1}^{p_{2,0}})|
=|\uparrow_{p_1}^{p_{2}}\uparrow_{p_1}^{p_{2,0}}|$ if
$|\uparrow_{p_1}^{p_{2}}\uparrow_{p_1}^{p_{2,0}}|=|p_2p_{2,0}|$. Then due to
Lemma \ref{3.12}, there is a natural 1-1 map (similar to the {\bf i} in Case 1)
$$\iota:M\to {\Bbb C\Bbb P}^{\frac n2}/\Bbb Z_2$$
such that $\iota([p_1p_2])=[\hat p_1\hat p_2]$ with $\uparrow_{\hat
p_1}^{\hat p_2}=i_{p_1*}(\uparrow_{p_1}^{p_2})$ for any $[p_1p_2]$
with $p_i\in M_i$.

\noindent{\bf Claim}:{\it\ $\iota$ is a continuous map (so it is a
homeomorphism).} Note that $\iota|_{M_i}$ is an isometrical
embedding for $i=1$ and 2, and that $\iota(M_i)=\hat M_i$. Moreover,
$\iota$ restricted to any convex spherical surface bounded by some
$[p_1p_2]$, $[p_1p_{2,0}]$ and $[p_2p_{2,0}]$ (here $p_i\in M_i$) is
an isometrical embedding (this is due to Lemma \ref{2.4}). It
is not hard to see that these together with (2.5) imply that $\iota$
is a continuous map.

From the above claim, we know that $\pi_1(M)\cong\Bbb Z_2$. Let
$\pi:\tilde M\to M$ be the Riemannian covering map. It suffices to
show that $\tilde M$ is isometric to ${\Bbb C\Bbb P}^{\frac n2}$.
Since $\iota$ is a homeomorphism and $\iota(M_i)=\hat M_i$,
we have that $\pi^{-1}(M_i)$ is connected (note that
$\hat M_i\ (={\Bbb C\Bbb P}^{m_i}/\Bbb
Z_2)\hookrightarrow {\Bbb C\Bbb P}^{\frac n2}/\Bbb Z_2$ induces an
isomorphism from $\pi_1(\hat M_i)$ to $\pi_1({\Bbb C\Bbb P}^{\frac
n2}/\Bbb Z_2)$),
and thus $$\pi^{-1}(M_i)\stackrel{\rm iso}{\cong}{\Bbb C\Bbb
P}^{m_i}$$ because $M_i\stackrel{\rm iso}{\cong}{\Bbb C\Bbb
P}^{m_i}/\Bbb Z_2$ and $M_i$ is totally geodesic in $M$
(which implies that $\pi^{-1}(M_i)$ is totally geodesic in
$\tilde M$). Moreover,
note that $|\tilde p_1\tilde p_2|\geq\frac\pi2$ for
any $\tilde p_i\in \pi^{-1}(M_i)$. Hence, $\tilde M$ satisfies
the conditions of the Main Theorem, so it follows from Case 1 that
$\tilde M$ is isometric to ${\Bbb C\Bbb P}^{\frac n2}$. \hfill$\Box$

%%%%%%%%%%%%%%%%%%%%%%%%%%% Appendix %%%%%%%%%%%%%%%%%%%%%%%%%%%%%%%%%%%%

\vskip8mm

\noindent{\bf \Large Appendix}

\vskip5mm

\noindent{\bf A.1.} On ${\Bbb C\Bbb P}^{m}/\Bbb Z_2$ ($m$ is odd)
with the canonical metric

How to get ${\Bbb C\Bbb P}^{m}/\Bbb Z_2$? We know that $$\Bbb
S^{2m+1}=\{(z_1,\cdots,z_{m+1})|z_i\in \Bbb C,\
|z_1|^2+\cdots+|z_{m+1}|^2=1\},$$ and $S^1$ can act on $\Bbb S^{2m+1}$
freely and isometrically (see (0.3)) through
$$S^1\times \Bbb
S^{2m+1}\to \Bbb S^{2m+1} \text{ defined by }
(e^{i\theta},(z_1,\cdots,z_{m+1}))\mapsto(e^{i\theta}z_1,\cdots,e^{i\theta}z_{m+1}).
$$
And if $m$ is odd, then there is an isometry of order $4$ on $\Bbb
S^{2m+1}$:
$$\varsigma:\Bbb
S^{2m+1}\to \Bbb S^{2m+1} \text{ defined by } (z_1,z_2,\cdots,
z_{2j-1},z_{2j},\cdots)\mapsto(-\bar z_2,\bar z_1,\cdots,-\bar
z_{2j},\bar z_{2j-1},\cdots).
$$
Note that $\varsigma$ induces a 2-order isometry $\hat\varsigma$
without fixed points on $\Bbb S^{2m+1}/S^1$. $(\Bbb
S^{2m+1}/S^1)/\langle\varsigma\rangle$ endowed with the induced
metric from the unit sphere $\Bbb S^{2m+1}$ is just the ${\Bbb C\Bbb
P}^{m}/\Bbb Z_2$ with the canonical metric. Moreover, for any odd
$m_i>0$ with $m_1+m_2=m-1$, ${\Bbb C\Bbb P}^{m_i}/\Bbb Z_2$ can be
embedded isometrically into ${\Bbb C\Bbb P}^{m}/\Bbb Z_2$ with
$$|q_1q_2|=\frac\pi2\text{ for any } q_i\in\Bbb{CP}^{m_i}/\Bbb Z_2.\eqno{\text{(A.1)}}$$
In other words, ${\Bbb C\Bbb P}^{m}/\Bbb Z_2$ is indeed  an example of the
Main Theorem.

As for even $m$, we have the following property.

\vskip2mm

\noindent{\bf Lemma A.1.} {\it $\Bbb Z_2$ can {\rm NOT} act on
${\Bbb C\Bbb P}^{m}$ freely by isometries when $m$ is even. }

\vskip2mm

\noindent{\bf Proof.} We will give the proof by the induction on $m$.
Obviously, when $m=0$, this is true (because $\Bbb{CP}^0$ contains
only one point). Now we assume that $m>0$, and that $\Bbb
Z_2\triangleq\langle\sigma\rangle$ acts on ${\Bbb C\Bbb P}^{m}$ by
isometries. Let $p$ be an arbitrary point in ${\Bbb C\Bbb P}^{m}$.
Note that $|p\sigma(p)|\leq\frac\pi2$, and that
$\sigma([p\sigma(p)])$ is also a minimal geodesic between $p$ and
$\sigma(p)$ for any $[p\sigma(p)]$. It follows that $\sigma$ fixes
the middle point of $[p\sigma(p)]$ if $|p\sigma(p)|<\frac\pi2$
(because there is a unique minimal geodesic between $p$ and
$\sigma(p)$ when $|p\sigma(p)|<\frac\pi2$). If
$|p\sigma(p)|=\frac\pi2$, then $\sigma$ preserves the set
$L\triangleq\{x\in \Bbb{CP}^m|x \text{ belongs to some }
[p\sigma(p)]\}$ which is isometric to $\Bbb{CP}^1$. Note that
$L^{=\frac\pi2}$ is isometric to
$\Bbb{CP}^{m-2}$ and totally geodesic (in $\Bbb{CP}^m$). Since
$\sigma$ preserves $L$, it has to preserves $L^{=\frac\pi2}$,
and thus $\sigma|_{L^{=\frac\pi2}}$  is an isometry. By the
inductive assumption, $\sigma$ has fixed points on
$L^{=\frac\pi2}$ (so on $\Bbb{CP}^m$). \hfill$\Box$

\vskip2mm

In fact, $\Bbb Z_2$ cannot act on ${\Bbb C\Bbb P}^{m}$ freely in the
sense of topology when $m$ is even ([Sa]).

\vskip4mm

\noindent{\bf A.2.} Proof of (3.4.1) in Lemma \ref{3.4} ([RW])

Since $\lambda_{p_1p_2}\equiv h$ for all $p_i\in M_i$, due to Lemma \ref{2.3}
it follows that, for the given
$p_1\in M_1$ and $p_2\in M_2$, there are $\varepsilon>0$ and a
neighborhood $V$ of $p_2$ in $M_2$ such that
$$\min\limits_{1\leq j\neq j'\leq h}
\{|(\uparrow_{p_1}^{p_2'})_j(\uparrow_{p_1}^{p_2'})_{j'}|
|(\uparrow_{p_1}^{p_2'})_j,(\uparrow_{p_1}^{p_2'})_{j'}\in
\Uparrow_{p_1}^{p_2'},\ p_2'\in V\}>\varepsilon. \eqno{\rm(A.2)}$$
Let $U=B(p_2,\frac{\varepsilon}{4})\cap V$. By Lemma
\ref{2.4}, for the given $[p_1p_2]$ and any $p_2'\in U$
$$\exists!\ [p_1p_2']\text{ such that }
|\uparrow_{p_1}^{p_2}\uparrow_{p_1}^{p_2'}|=|p_2p_2'|.
\eqno{\rm(A.3)}$$ Note that we need only to prove that
$|\uparrow_{p_1}^{p_2^1}\uparrow_{p_1}^{p_2^2}|=|p_2^1p_2^2|$ for
all $p_2^1,p_2^2\in U$, where $\uparrow_{p_1}^{p_2^j}$ is the
direction of the $[p_1p_2^j]$ found in (A.3). If this is not true,
by Lemma \ref{2.3} there is another minimal geodesic $[p_1p_2^2]'$
between $p_1$ and $p_2^2$ such that
$|\uparrow_{p_1}^{p_2^1}(\uparrow_{p_1}^{p_2^2})'|=|p_2^1p_2^2|$.
However,
$|\uparrow_{p_1}^{p_2^2}(\uparrow_{p_1}^{p_2^2})'|\leq|p_2^2p_2|+
|p_2p_2^1|+|p_2^1p_2^2|<\varepsilon$, which contradicts (A.2).
\hfill$\Box$

\vskip4mm

\noindent{\bf A.3.} On the metric of $U_1*U_2$ in (3.4.2)

Let $X$ and $Y$ be two Alexandrov spaces with curvature $\geq 1$
(especially two Riemannian manifolds with sectional curvature $\geq
1$). The canonical metric on the join space $X*Y\triangleq X\times
Y\times[0,\frac{\pi}{2}]/\sim$, where
$(x,y,t)\sim(x',y',t')\Leftrightarrow t=t'=0$ and $x=x'$ or
$t=t'=\frac{\pi}{2}$ and $y=y'$, is defined as follows ([BGP]):
$$\cos|p_1p_2|=\cos t_1\cos t_2\cos|x_1x_2|+\sin t_1\sin t_2\cos|y_1y_2|$$
(with $|p_1p_2|\leq\pi$) for any $p_i\triangleq[(x_i,y_i,t_i)]\in X*Y$. It can be proved
([BGP]) that $X*Y$ (endowed with such a metric) is also an
Alexandrov space with curvature $\geq 1$ and $\dim(X*Y)=\dim(X)+\dim(Y)+1$ (especially, it is the unit
sphere if both $X$ and $Y$ are unit spheres).

\vskip4mm

\noindent{\bf A.4.} On the case where $n_1,n_2>0$ in the Main Theorem

\vskip2mm

\noindent{\bf Proposition A.4.} {\it In the Main Theorem,
if $n_1,n_2>0$, then either $M\stackrel{\rm
iso}{\cong}\Bbb S^n$ or $\Bbb{RP}^n$, or $M_1=M_2^{\geq\frac\pi2}\ (=M_2^{=\frac\pi2})$
and $M_2=M_1^{\geq\frac\pi2}\ (=M_1^{=\frac\pi2})$.}

\vskip2mm

\noindent{\bf Proof.} Let $\bar M_1$ denote $M_2^{\geq\frac\pi2}$.
It suffices to show that $M\stackrel{\rm
iso}{\cong}\Bbb S^n$ or $\Bbb{RP}^n$ if $M_1\neq \bar M_1$.
Note that $\bar M_1=M_2^{=\frac\pi2}$ by Lemma \ref{2.1}.
Then according to [RW] (cf. Footnote 5), $\dim(\bar M_1)=n_1+1$, and
$\lambda_{q_1q_2}\equiv \bar h<+\infty$ for
all $q_1\in \bar M_1^\circ$ and  $q_2\in M_2$ (note that
$\bar M_1$ may have nonempty boundary). Fix an arbitrary $[q_1q_2]$ with
$q_1\in \bar M_1^\circ$ and  $q_2\in  M_2$, and consider
$f_{[q_1q_2]}:M_2\rightarrow (\Sigma_{q_1}\bar M_1)^{=\frac\pi2}\ (=\Bbb S^{n_2})$. By (3.4.1), we
conclude that $\sec_{M_2}\equiv 1$ or $n_2=1$, so
$\lambda_{p_1p_2}\equiv h<+\infty$ for
all $p_i\in M_i$ by Proposition \ref{3.1}. Hence,
it follows from Lemma \ref{3.3} and its proof that
$M\stackrel{\rm iso}{\cong}\Bbb S^n$ or $\Bbb{RP}^n$.
\hfill$\Box$

%%%%%%%%%%%%%%%%%%%%%%%%%%%%%%%%%%%%%%%%%%%%%%%%%%%%%%%%%%%%%%%

\noindent School of Mathematical Sciences (and Lab. math. Com.
Sys.), Beijing Normal University, Beijing, 100875
P.R.C.\\
e-mail: suxiaole$@$bnu.edu.cn; wyusheng$@$bnu.edu.cn

\vskip2mm

\noindent Mathematics Department, Capital Normal University,
Beijing, 100037 P.R.C.\\
e-mail: 5598@cnu.edu.cn

\end{document}